\documentclass{amsart}

\usepackage{amssymb}
\usepackage{amscd}
\usepackage{latexsym}
\usepackage{amsfonts}
\usepackage{amsmath}
\usepackage{graphicx}
\usepackage{texdraw}
\newtheorem{prop}{Proposition}[section]

\newtheorem{remark}{Remark}[section]
\newtheorem{lemma}{Lemma}[section]

\newtheorem{theorem}{Theorem}[section]

\newtheorem{definition}{Definition}[section]

\newtheorem{claim}{Claim}[section] 

\newtheorem{conjecture}{Conjecture}[section]

\newcommand{\im}{{\rm  Im}}

\newcommand{\spane}{{\rm  Span}}

\newcommand{\tr}{{\rm  Tr}}

\newcommand{\Hom}{{\rm  Hom} }

\newcommand{\expo}{\exp}

\newcommand{\rank}{{\rm   Rank}}
\newcommand{\ra}{\rightarrow}

\newcommand{\Diff}{{\rm  Diff} }
\newcommand{\Ad}{{\rm   Ad} }

\newcommand{\Det}{{\rm   Det}}

\newcommand{\ch}{{\rm   ch} }

\newcommand{\id}{{\rm   Id}}

\newcommand{\Aut}{{{\rm   Aut} }}

\newcommand{\graph}{{\rm   graph}}
\newcommand{\sign}{{\rm   sign} }
\newcommand{\alg}{{\rm   alg}}

\newcommand{\eq}{{\rm   eq} }
\newcommand{\Td}{{\rm   Td}}
\newcommand{\Id}{{\rm   Id} }
\newcommand{\CS}{{\rm   CS} }
\newcommand{\Ch}{ {\rm Ch}}
\newcommand{\Pic}{{\rm   Pic} }
\newcommand{\R}{  {\mathbb R}} 
\newcommand{\C}{{\mathbb C}} 
 
\newcommand{\Z}{{\mathcal Z}}
\newcommand{\f}{ {\tilde f}}

\newcommand{\M}{{\mathcal M}}
\newcommand{\Q}{ {\mathbb Q}}

\begin{document} 


\title[Quantum invariants of finite order mapping tori I.]{The Witten-Reshetikhin-Turaev invariants of finite order mapping tori I.}

\author{J\o rgen Ellegaard Andersen}
\address{Center for Quantum Geometry of Moduli Spaces\
        University of Aarhus\\
        DK-8000, Denmark}
\email{andersen@imf.au.dk}

\thanks{This paper is an update version of our preprint "The Witten Invariant of Finite Order Mapping Tori I.", Aarhus University, Department of Mathematical Sciences, Preprint series 1995, No. 13.} 

\thanks{Supported in part by NSF grant
DMS-93-09653, while the author was visiting University of California, Berkeley. }

\thanks{Supported in part by the center of excellence grant "Center for quantum geometry of Moduli Spaces" from the Danish National Research Foundation. }

\begin{abstract}

{\noindent We formulate the Asymptotic Expansion Conjecture for the Witten-Reshetikhin-Turaev quantum invariants of closed oriented three manifolds. For finite order mapping tori, we study these quantum invariants via the geometric gauge theory approach to the corresponding quantum representations and prove, using a version
 of the Lefschetz-Riemann-Roch Theorem due to Baum, Fulton, MacPherson \& Quart,
that the quantum
invariants can be expressed as a sum over the components of the moduli space of flat 
connections on the mapping torus. 
Moreover, we
show that the term corresponding to a component is a polynomial in the level
$k$, weighted by a complex phase, which is $k$ times the Chern-Simons invariant
corresponding to the component. We express the coefficients of these polynomials
in terms of cohomological pairings on the fixed point set of the moduli space of
flat connections on the surface.  We explicitly describe the fixed point set in
terms of moduli spaces of the quotient orbifold Riemann surface and for the
smooth components we express
the aforementioned  coefficients in terms of the known generators of the
cohomology ring. We provide an explicit formula in terms of the Seifert invariants of the mapping torus for the contributions from each of the smooth components.
We further establish that the Asymptotic Expansion Conjecture  and  
the Growth Rate Conjecture for these finite order mapping tori.}

\end{abstract}

\maketitle

\section{Introduction}

In this paper we shall describe how the Witten-Reshetikhin-Turaev quantum invariants of an interesting class of 3-manifolds can be computed using the rigorous geometric approach. By $\Z^{(k)}_G$ we will denote the Reshetikhin-Turaev TQFT at level $k$ for the semisimple, simply connected Lie group $G$. For the construction of this TQFT please see \cite{RT1}, \cite{RT2}, \cite{T}, \cite{TW} and \cite{KL}. 

From general path integral techniques (see e.g. \cite{Witten}, \cite{FG}, \cite{J1}, \cite{Ro1}, \cite{Ro2}, \cite{LR}, \cite{LZ}, \cite{M}) we expect that the Witten-Reshetikhin-Turaev invariants will have a asymptotic expansion in the level of the theory.  There are a number of formulations of this expectation in the literature, which expresses this as a sum over flat $G$-connections on the three manifold. These formulations come naturally from  path integrals, however this expresses the asymptotics in terms of sum of terms which can not be uniquely recovered from the asymptotic expansion, e.g. sums over flat $G$-connections. Let us here state a very precise version of this Conjecture, which also appeared in \cite{AAEC} for a general oriented compact three manifold $X$. It has the property that each term in the expansion can be uniquely recovered from the quantum invariant $\Z^{(k)}_G(X)$.

\begin{conjecture}[Asymptotic Expansion Conjecture
(AEC)]\label{conj:AEC}
There exist constants (depending on X) $d_j \in \frac12 {\mathbb Z}$ and
$b_j \in \C$ for $j=0,1, \ldots, n$ and $a_j^l \in  \C$ for
$j=0,1, \ldots, n$,\, $l=1,2,\ldots$ such that the asymptotic
expansion of $\Z^{(k)}_G(X)$ in the limit $k\ra \infty$ is
given by
\[\Z^{(k)}_G(X) \sim \sum_{j=0}^n e^{2\pi i k q_j} k^{d_j}
b_j \left( 1 + \sum_{l=1}^\infty a_j^l k^{-l/2}\right), \]
where $q_0 = 0, q_1, \ldots q_n$ are the finitely many
different values of the Chern--Simons functional on the space
of flat $G$--connections on $X$.
\end{conjecture}

Here {\bf $\sim$} means {\bf asymptotic expansion} in the
Poincar\'{e} sense, which means the following: Let
\[d = \max\{d_0,\ldots,d_n\}.\]
Then for any non-negative integer $L$, there is a $c_L \in \R$
such that
\[\left| \Z^{(k)}_G(X) - \sum_{j=0}^n e^{2\pi i k q_j} k^{d_j}
b_j \left( 1 + \sum_{l=0}^L a_j^l
k^{-l/2}\right) \right| \leq c_L k^{d-(L+1)/2}\]
for all levels $k$. Of course such a condition only puts
limits on the large $k$ behavior of $ \Z^{(k)}_G(X)$.

This expansion tends to organize more beautiful if one expands in the parameter $r=k+h$, where $h$ is the dual Coexeter number of the group. This is discussed in a number of the references provided below. We comment further on this below.

A little simple argument shows that there exists
at most one list of numbers $n \in \{0,1,\ldots\}$,
$q_0,q_1,\ldots,q_n \in \R \cap [0,1[$, $d_j \in \frac12 {\mathbb Z}$
and $b_j \in \C$ for $j=0,1,\ldots, n$
and $a_j^l \in  \C$ for $j=0,1, \ldots, n$,\,
$l=1,2,\ldots$ such that $0=q_0 < q_1 < \ldots < q_n$ and such
that the large $k$ asymptotic expansion of $\Z^{(k)}_G(X)$ is given as above.
This
implies that if the quantum invariant $\Z^{(k)}_G(X)$ has an
asymptotic expansion of the above form, then
$q_j$'s, $d_j$'s, $b_j$'s and $a_j^l$'s are all uniquely
determined by $\Z^{(k)}_G(X)$,
hence they are also topological invariants of $X$. As stated above
the AEC already includes the claim that the $q_j$'s are the
Chern--Simons invariants. There are also conjectured topological
formulae for the $d_j$'s, and $b_j$'s (see e.g. the above given references). Let us here just recall the
conjectured formula for the $d_{j}$'s.

For a flat $G$--connection $A$ on the $3$--manifold $X$,
denote by $h_{A}^{i}$ the dimension of
the $i'$th $A$-twisted cohomology groups of $X$ with Lie algebra coefficients.

\begin{conjecture}[Growth Rate Conjecture]\label{conj:d}
Let $\M^{(j)}$ be the union of components of the moduli space
of flat $G$--connections on $X$ which has Chern--Simons
value $q_{j}$. Then
$$
d_{j} = \frac{1}{2} \max_{A \in \M^{(j)}}
            \left( h_{A}^{1} - h_{A}^{0} \right),
$$
where $\max$ here means the maximum value
$h_{A}^{1} - h_{A}^{0}$ attains on a Zariski open subset
of $\M^{(j)}$.
\end{conjecture}

One expects that there are expressions for each of the
$a_j^l$'s in terms of sums over Feynman diagrams of certain
contributions determined by the Feynman rules of the
Chern--Simons theory. This has not yet been worked out in
general, except in the case of an acyclic flat connection and
the case of a smooth non-degenerate component of the moduli
space of flat connections by Axelrod and Singer,
cf.\ \cite{[4]}, \cite{[5]}, \cite{[3]}. See also the work of Marino \cite{M} and references therein. In the end of this introduction we will provide some references to results, which establishes the above two Conjectures in various cases.

The AEC, Conjecture \ref{conj:AEC}, however offers in a sense a
converse point of view, where one seeks to derive the final
output of perturbation theory after all cancellations have
been made (i.e.\ collect all terms with the same Chern--Simons
value and same power of $k$). This seems actually rather reasonable in this case,
since the exact invariant is known explicitly.

In this paper we consider the three manifolds which are mapping tori. Let us now recall a formula for the quantum invariants of such. 

The Witten-Reshetikhin-Turaev TQFT provide representations of central extensions of mapping class groups (see also \cite{MRob}): Suppose $\Sigma$ is a closed oriented surface of genus $g$ and let $\Gamma$ be the mapping class group of  $\Sigma$. Choose a Lagrangian subspace $L$ of $H_1(\Sigma,\R)$. We denote by $\Gamma_r$ the rigged mapping class group of $(\Sigma,L)$, as constructed in \cite{T} and \cite{Walker} (see also Section 2). Recall that $\Gamma_r$ is a ${\mathbb Z}$-central extension of $\Gamma$. The TQFT $\Z^{(k)}_G$ assigns a vector space to the pair $(\Sigma,L)$, which we simply denote $\Z^{(k)}_G(\Sigma)$, and a representation
$$ \Z^{(k)}_G : \Gamma_r \rightarrow \Aut(\Z^{(k)}_G(\Sigma)).$$

Suppose $f$ is an orientation preserving diffeomorphism of a closed surface
\(\Sigma\). We can then form the mapping torus of $f$: 
$$\Sigma_{f} = (\Sigma\times {\mathbb R})/{\mathbb Z},$$
where  $p\in {\mathbb Z}$  acts by $p(x,t) = (f^{-p}(x),t+p)$ for all $(x,t)\in \Sigma\times {\mathbb R}$. The product orientation on $\Sigma\times {\mathbb R}$ induces an orientation on $\Sigma_{f}$. We will 
review how the choice of a 2-framing of $\Sigma_{f}$ gives an
element in the central extension of the mapping class group. Hence, if we use
the Atiyah 2-framing of $\Sigma_f$, we get an element $\tilde{f}$ of this central
extension. This will be explained in details in Section \ref{Z}
and \ref{DEF}. The
axioms for a TQFT imply that the Witten-Reshetikhin-Turaev quantum invariant \(\Z^{(k)}_G(\Sigma_{f})\) is given by:

\[\Z^{(k)}_G(\Sigma_{f}) = \tr (\Z^{(k)}_G(\f) : \Z^{(k)}_G(\Sigma) \rightarrow \Z^{(k)}_G(\Sigma)).\] 
We  thus see that the character contains
the information about the Witten-Reshetikhin-Turaev quantum invariants  of an interesting class of $3$-manifolds.

In this paper we will use the gauge theory construction of the quantum representations $Z^{(k)}_G$ of the mapping class group as we did in \cite{A1,A2}. 
We establish in the series of papers \cite{AU1, AU2, AU3} together with \cite{AU4} in preparation joint with Kenji Ueno, that the gauge theory construction of the 
quantum invariants for $G=SU(n)$ coincides with the ones first constructed by Reshetikhin-Turaev in \cite{RT1, RT2, T, TW} or equivalently by Blanchet-Habegger-Masbaum-Vogel in \cite{BHMV1, BHMV2, B}. We stress that the results presented in this paper is however independent of this result, since we work entirely on the gauge theory side and only with the gauge theory definition of the quantum representations.

Let us recall the gauge theory construction of the vector space $Z^{(k)}_G(\Sigma)$.
One applies geometric quantization to the moduli space of flat $G$-connections on the surface $\Sigma$, so as to obtain a vector bundle  \(Z\)  (we suppress the dependence of the Lie group $G$ and the level $k$) with a flat connection over Teichm\"{u}ller
space $T$.  We denote the space of  covariant constant sections
of this bundle by $Z^{(k)}_G(\Sigma)$.  Moreover $\Gamma_r$
acts on \(Z\) and preserves the connection, so we get a representation
$$Z^{(k)}_G : \Gamma_r \ra \Aut(Z^{(k)}_G(\Sigma)).$$
We observe that if $f$ is finite order in the mapping class group, then it preserves a point
\(\sigma\) in Teichm\"{u}ller space and
\[Z^{(k)}_G(\Sigma_f) := \tr Z^{(k)}_G(\tilde{f}) = \tr (\f : Z_\sigma \rightarrow Z_\sigma).\]

In Section \ref{Z} we shall explain in detail
how $Z$ is constructed. The outline of the construction is as follows.
Pick a compact simple simply-connected Lie group $G$ with a suitable
normalized biinvariant inner product. Let ${\mathcal M}$ be the moduli space of flat $G$ connections on $\Sigma$. This is a stratified symplectic space on which the mapping class group acts. For a $\sigma\in T$, consider the moduli space ${\mathcal
M}_\sigma$ of semi-stable $G^{\Bbb C}$ bundles over $\Sigma_\sigma$. One can identify ${\mathcal M}$ and ${\mathcal M}_\sigma$ as stratified spaces, but the moduli
space ${\mathcal M}_\sigma$ further has the structure of a normal projective variety.
There is a natural construction of the determinant line bundle ${\mathcal L}_\sigma$ over ${\mathcal
M}_\sigma$ whose first Chern class is represented by the K\"{a}hler form on
${\mathcal M}_\sigma$.  We consider the vector bundle $\hat{Z}$ over $T$, whose fiber over $\sigma\in T$ is
$$\hat{Z}_\sigma = H^0({\mathcal
M}_\sigma,
{\mathcal L}_\sigma^k).$$

As will be explained in Section \ref{Z}, there is a natural lift of the action of the mapping class group $\Gamma$ to the family of line bundles ${\mathcal L}_\sigma$,  $\sigma\in T$.
Thus we see there is a natural action of $\Gamma$ on $\hat{Z}$.
By the results of S. Axelrod, S. Della Pietra \& Witten \cite{ADW}, N. Hitchin \cite{Hitchin2} and G. Faltings \cite{Faltings} the bundle $\hat{Z}$ has a natural  mapping class group invariant connection, which is projectively flat.
Since this connection is not flat, this bundle needs to be modified by a line bundle with a suitable connection.
Let ${\mathcal L}_D$ be the determinant 
bundle over Teichm\"{u}ller space. Let $\zeta\in \Q$ be the central change for $G$ (see formula 
(\ref{cc}) below). 
As it is described in \cite{Walker} on can construct any fractional power of ${\mathcal L}_D$ 
over $T$ and the rigged mapping class group
$\Gamma_r$ naturally acts on it. 
By definition, the fiber over $\sigma$ of the above mentioned vector bundle $Z$, is
\begin{equation}
Z_\sigma = H^0({\mathcal M}_\sigma,
{\mathcal L}_\sigma^k)\otimes {\mathcal L}_{D,\sigma}^{-\frac{1}{2}\zeta}.\label{F1}
\end{equation}
As it was discussed in \cite{AU2}, $\Gamma_r$ acts naturally on the line bundle ${\mathcal L}_{D}^{-\frac{1}{2}\zeta}.$

We will throughout the rest of the introduction on assume that our mapping class $f$ is finite order and that $\tilde{f}$ is a lift of it to $\Gamma_r$. Then $f$ preserves a point $\sigma$ in $T$ and we get that
\begin{equation}
 \tr Z^{(k)}_G(\tilde{f}) = \tr(f: H^0({\mathcal M}_\sigma,{\mathcal L}_\sigma^k) \rightarrow H^0({\mathcal
M}_\sigma,{\mathcal L}_\sigma^k))\tr(\tilde{f} : {\mathcal
L}_{D,\sigma}^{-\frac{1}{2}\zeta} \rightarrow {\mathcal
L}_{D,\sigma}^{-\frac{1}{2}\zeta}).\label{F34}
\end{equation}
Let us introduce the following notation
$$ \tau^{(k)}_G(f) = \tr(f: H^0({\mathcal M}_\sigma,{\mathcal L}_\sigma^k) \rightarrow H^0({\mathcal
M}_\sigma,{\mathcal L}_\sigma^k))$$
and 
$$ \Det(f)^{-\frac{1}{2}\zeta} = \tr(\tilde{f} : {\mathcal
L}_{D,\sigma}^{-\frac{1}{2}\zeta} \rightarrow {\mathcal
L}_{D,\sigma}^{-\frac{1}{2}\zeta}),$$
thus we have that

$$ Z^{(k)}_G(\Sigma_f) =  \Det(f)^{-\frac{1}{2}\zeta}  \tau_G^{(k)}(f).$$

The {\em framing correction} $\Det(f)^{-\frac{1}{2}\zeta}$ is calculated in Theorem \ref{fct} in Section
\ref{FC}. Formula (\ref{FCT}) below expresses the framing correction term as a function of $m$ and the Seifert invariants of $\Sigma_f$ as an explicit power series $1/k$. It is remarkable that the Atiyah 2-framing choice is precisely such that
the expression for this term simplifies substantially. 

It is  the factor $ \tau_G^{(k)}(f)$ which is by far the most interesting. 
By
general theory one can show that the higher cohomology groups of the line
bundle ${\mathcal L}_\sigma^k$ vanish. This means we are exactly in the situation, 
where we can apply
the Lefschetz-Riemann-Roch Theorem for finite automorphism of projective
varieties due to Baum, Fulton, MacPherson \& Quart. That means we get an expression
 for $\tau_G^{(k)}(f)$, as a sum over the components of the fixed point set
of the automorphism induced by $f$ on ${\mathcal M}_\sigma$. (See Section
\ref{FPS}.)

The upshot of this calculation is a formula for $\tau_G^{(k)}(f)$  of any finite order diffeomorphism\footnote{Recall that a Seifert fibred
manifold is  a mapping torus of a finite order diffeomorphism iff it has zero
S.F.-Euler number, hence we get a large class of the Seifert fibred spaces from finite order elements of the mapping class group. See Section \ref{SI}.}, which is expressed 
as a sum over the components of the
fixed point set $|{\mathcal M}_\sigma|$ of $f$ on ${\mathcal M}_\sigma$. Let us denote the set of components $C$. For each $c\in C$, we let $|{\mathcal M}_\sigma|^c$ be the component 
subset of $|{\mathcal M}_\sigma|$. We will see that the contribution from $|{\mathcal M}_\sigma|^{c}$
is a root of unity to the power $k$ multiplied by a {\em polynomial} in k, say $P_c$, where
each of the coefficients of $P_c$ are expressed as certain cohomological pairings on the
fixed point set (see Section \ref{L}). The degree of $P_c$ we denote $d_c$.

By
understanding the relation between the moduli space of flat connections on the
mapping torus and the fixed point set of $f$ on ${\mathcal M}_\sigma$ (see Section
\ref{3MF}),
we can identify the contribution 
from the different components of the moduli space of flat connection on the
mapping torus. In fact, suppose ${\mathcal M}(\Sigma_f)$ is the moduli space of flat $G$-connections. We get a 
map from ${\mathcal M}(\Sigma_f)$ to $|{\mathcal M}_\sigma|$ by restriction to $\Sigma\times \{0\}$, which of course takes components
to components. Let ${\mathcal M}(\Sigma_f)^c$ be the union of components of ${\mathcal M}(\Sigma_f)$ 
which maps to $|{\mathcal M}_\sigma|^c$.

The term in the sum corresponding to $c$ is given as the exponential of $2\pi i k$
times the Chern-Simons invariant $\CS(\Sigma_f,c)$ corresponding to the collection of components ${\mathcal M}(\Sigma_f)^c$ (we will see by Lemma \ref{TCS} that the Chern-Simons invariant only depends on $c$) multiplied by $P_c(k)$. 

We thus see that every component of the moduli space
of flat connections on $\Sigma_f$ only contributes finitely many loop corrections to the invariant, 
modulo of course the framing term $\Det(f)^{-\frac{1}{2}\zeta}$.
The framing term gives an over all common factor, which is a power series $1/k$ given explicitly in formula (\ref{FCT}) below. 

See Theorem \ref{TWF} for a detailed statement of the result. Let us here summarize our main result in the following Theorem

\begin{theorem}
The AEC holds for all finite order mapping tori. In particular for a finite order mapping class $f$, 
the Witten-Reshetikhin-Turaev invariant of the mapping torus
$\Sigma_f$ is given as the following finite sum 
\[Z^{(k)}_G(\Sigma_f) = \Det(f)^{-\frac{1}{2}\zeta}
\sum_{c\in C} \expo(2 \pi i k\CS(\Sigma_f,c))k^{d_c} P_c(k^{-1}),\] 
where $ \Det(f)^{-\frac{1}{2}\zeta}$ is the power series in $1/k$ given by formula (\ref{FCT}) below and the polynomial $P_c$ is of order $d_c$ and each of its coefficients are expressed as certain intersection numbers on the component $|{\mathcal M}_\sigma|^c$ of the fixed point set $|{\mathcal M}_\sigma|$ as stated in Theorem \ref{TWF} below. 
\end{theorem}

By further analyzing dimensions of twisted cohomology groups (see Lemma \ref{dc}), we get the following Theorem

\begin{theorem}
For each of the connected components $c\in C$ we get that
$$
d_c = \frac{1}{2} \max_{A \in \M(\Sigma_f)^c}
            \left( h_{A}^{1} - h_{A}^{0} \right),
$$
where $\max$ here means the maximum value
$h_{A}^{1} - h_{A}^{0}$ attaines on a Zariski open subset
of $\M(\Sigma_f)^c$, hence the Growth Rate Conjecture holds for finite order mapping tori.
\end{theorem}

\begin{remark}
The fact that one really needs the generic maximum in formula (\ref{dF}) is
illustrated by the following simple example. If one considers the mapping torus
of $-\id : T^2 \rightarrow T^2$, then it is easily seen that $Z(T^2_{-\id})=
Z(T^2\times S^1) = k+1$, however an explicit calculation in this particular 
case shows that
the moduli space of flat $SU(2)$ connections on
$T^2_{-\id}$ is connected and the maximal value of $\frac{1}{2}(h_{A}^{1} - h_{A}^{0} )$ is $2$, but  that the
generic maximum is $1$ as it should be.
\end{remark}

For each smooth component we express the coefficient in terms
of known generators of the cohomology ring of that component of the moduli space
of flat connections on the mapping torus. The result is precisely stated in
Theorem \ref{SAET2}. From the explicit expression obtained, we observe that the
contribution from the smooth components of the moduli space only depends on $m$ and the
Seifert invariants of the mapping torus.

Throughout this paper, an attempt has been made to emphasize the use of flat
bundles as opposed to semi-stable parabolic bundles in order to keep the
connection with flat bundles on the mapping torus more transparent.

The asymptotic expansion of the Witten-Reshetikhin-Turaev quantum invariants has been studied by a number of authors and various results has been obtained for certain classes of closed three manifolds. Let us here comment on the works which are closest in relation to our work.

Freed and Gompf did some computer calculations of the asymptotics of the Witten-Reshetikhin-Turaev-invariants for lens spaces and certain Seifert fibered spaces in \cite{FG} providing evidence for a relation between the quantum invariants and classical Chern-Simons gauge theory invariants, which are proved in this paper for finite order mapping tori of surfaces of genus at least $2$.

In \cite{G} Garoufalidis considered lens spaces and verified the AEC for these. He also did some calculations for certain Seifert fibered homology spheres.
In \cite{J1} Jeffrey also consider the case of lens spaces and gave a proof of the AEC for these. In the same paper she further considered certain (not including the finite order ones) mapping tori of genus one surfaces and she verified the leading order asymptotic behavior of the quantum invariants for these. Rozansky studied in \cite{Ro1} the asymptotic expansion of the quantum invariants for Seifert fibered manifolds. However his approach in that paper did not apply to the Seifert fibered manifolds considered in this paper, since he needed to assume that the orbifold Euler characteristic was non-vanishing, which is not satisfied for mapping tori of finite order. The same remark applied to his subsequent work joint with Ruth Lawrence \cite{LR} and her work joint with D. Zagier \cite{LZ} dedicated to the study of quantum invariants. Rozansky did study the asymptotic expansion of the quantum invariants for $SU(2)$ for all Seifert fibered spaces over orientable surfaces in \cite{Ro2}. Rozansky starts with the combinatorial formulae for the Witten-Reshetikhin-Turaev invariants for these manifolds, but it is unclear to us to what extend the calculations presented in that paper are entirely rigorous. Furthermore, that paper only claims to get an asymptotic expansion in the level of the quantum invariants. Our results are exact results with no asymptotic approximation needed, since we show that the expansion is finite (modulo the framing correction term, which is provided explicitly). However, in Section 5 of that paper, Rozansky does state formulae for contributions to the asymptotic expansions of the quantum invariants, which he attributes to certain connected smooth components of the moduli space of flat $SU(2)$-connections. He further states expressions for these contributions in terms of intersection pairings on these moduli spaces. Hence his formulae shares a formal resemblance to ours. We have however not been able to reconcile our formulae with his completely and it does seem non-trivial to match these formulae up. Nevertheless if we consider Rozansky's expressions (5.20) (H=0)  in Proposition 5.3 in \cite{Ro2} it shares some formal resemblance to our formula in Theorem \ref{SAET2} in this paper, when we keep in mind that the $\Td$-class is the $\hat A$-class multiplied by $\exp(\frac12 c_1)$, which of course is related to the level shift $k\mapsto k+ h$ as mentioned above. It would indeed be interesting to establish the equivalence of the two expressions for smooth components.

In \cite{M}, Marino relates the asymptotics of quantum invariants of Homology sphere Seifert fibered spaces to asymptotics of matrix models. As of now, there is no direct overlap between Marino's work and this paper since none of the mapping tori are homology spheres. It would however be very interesting to understand if Marino's work can be generalized to our setting.

In \cite{Kold1} Hansen  investigated the asymptotic expansion of all Seifert fibered spaces and the existence of such was proved for this class of manifolds in \cite{Kold4}. This work was further extended in \cite{Kold2} and \cite{Kold3} joint with Takata. In this paper, we give however a formula for the full expansion in terms of cohomology pairings on the moduli space, something which is not covered by the works of Hansen nor his work joint with Takata. 

In \cite{BW} Beasley and Witten considered the path integral formula for these quantum invariants. Since they are working with path integrals, their work is per se not rigorous, however they provide path integral arguments for the fact that the perturbation expansion of these invariants are finite (modulo the framing correction term). We provide a mathematical proof of that in this paper.

In \cite{AH}, the author has jointly with B. Himpel applied the results of this paper to identify the leading order in the asymptotics of the quantum invariants entirely in terms of classical topological invariants of the mapping torus. That paper establishes, as an application of the results of this paper, that the leading order asymptotics of the quantum invariants of finite order mapping tori 
is given as predicted by path integral considerations originally due to Witten \cite{Witten}, and elaborated upon in \cite{FG}, \cite{J1}, \cite{Ro1}, \cite{Ro2}, \cite{LR}, \cite{LZ}, \cite{M}.  

The author would like to thank Simon K. Donaldson, Gregor Masbaum, Michael Thaddeus
and Kevin Walker for very useful discussion on this project.
Most of this work was done at the Department of Mathematics, University of
California, Berkeley and the author would like to thank the Department for its
hospitality.

\section{The gauge theory construction of the functor $Z^{(k)}_G$.}

\label{Z} 

We shall in this Section briefly review the geometric gauge theory construction of a
functor $Z^{(k)}_G$ from the category of closed rigged surfaces of genus $g>1$ to the
category of finite dimensional complex vector spaces. This construction is due
to S. Axelrod, S. Della Pietra \& E. Witten; G. Segal; M. Atiyah; N. Hitchin;
A.A. Beilinson \& D. Kazhdan; A.A. Beilinson \& V.V. Schechtman; A. Tsuchyia, 
K. Ueno \& Y. Yamada;
S. E. Cappell, R. Lee \& E. Y. Miller;
G. Faltings (See
\cite{Hitchin2}, \cite{ADW}, \cite{Segal}, \cite{BK}, \cite{BS}, \cite{TUY},
 \cite{CML} and \cite{Faltings}).

As briefly outlined in the introduction, the functor $Z^{(k)}_G$ is constructed by
applying the machinery of geometric K\"{a}hler quantization to the moduli space
of flat connections in some principal bundle over the surface. So we need to 
construct a Hermitian line bundle with connection over this moduli space, and a
lift of the action of the mapping class group to this bundle. This construction
for $SU(2)$ is due to Ramadas, Weitsman and Singer (see \cite{RSW}). We refer
to D. Freed (see \cite{F}) for the construction in the general case. Let us
here just recall the choices involved.

Let $G$ be a simple and simply-connected compact Lie group. Under these
assumptions any principal $G$-bundle over a compact surface is trivial.
Choose an invariant inner product\footnote{Following Freed (see \cite{F}), we
normalize this invariant  inner product by requiring that the closed 3-form
$-\frac{1}{6}\langle\theta\wedge[\theta\wedge\theta]\rangle$ represents a primitive
integral class in $H^3(G,{\Bbb R}),$ where $\theta$ is the Maurer-Cartan form.} $\langle\cdot,\cdot\rangle$ on
the Lie algebra of $G$. 
Let $\Sigma$ be an oriented closed surface. Let ${\mathcal A}_P$ be the space 
of connections in the trivial principal $G$-bundle $P$ over $\Sigma$.
Let ${\mathcal M}$ be the moduli space of flat $G$-connections in $P$.
If ${\mathcal M}'$ denotes the subset of ${\mathcal M}$ corresponding to
equivalence classes of irreducible connections, then ${\mathcal M}'$ is a smooth
symplectic manifold with the symplectic structure induced from that on
${\mathcal A}_P$
by symplectic reduction.  
We have a natural action of $\Gamma$ on ${\mathcal M}$ induced by
pull back. 
In this situation Freed describes in \cite{F}, how one constructs a smooth Hermitian
line bundle with a connection ${\mathcal L}'$ over $\M'$. He further describes how the construction can be extended to produce a topological line bundle ${\mathcal L}$ over ${\mathcal M}$, with a natural lift of the 
$\Gamma$-action. One immediately observes that this natural action has the property that any mapping class acts by the identity over the trivial flat connection.

We now need a K\"{a}hler polarization on ${\mathcal M}$. Suppose we are given
a complex structure on $\Sigma$, i.e. a point $\sigma$ in Teichm\"{u}ller space
$T$ of $\Sigma$. By a Theorem of Ramanathan ${\mathcal M}$ is isomorphic to the moduli space ${\mathcal M}_\sigma$ of
semi-stable holomorphic $G^{\Bbb C}$-bundles over $\Sigma_\sigma$, and as such has the
structure of a complex projective algebraic variety (See Theorem 7.1 in
\cite{R}. This Theorem was first proved in the $SU(N)$-case by Narashimhan and
Seshadri. See \cite{NS}.).

Using the Hermitian connection in ${\mathcal L}'$, we get a
unique holomorphic structure on ${\mathcal L}'$ over the complex
manifold ${\mathcal M}'$. We denote the resulting holomorphic line bundle ${\mathcal L}'_\sigma$. In \cite{DN}, Drezet and Narasimhan
 proved that the
restriction map 
$$|_{{\mathcal M}'_\sigma}: \Pic({\mathcal M}_\sigma) \rightarrow \Pic({\mathcal M}_\sigma')$$
 is an isomorphism, hence there
is a unique holomorphic line bundle ${\mathcal L}_\sigma$ over 
${\mathcal M}_\sigma$ such that ${\mathcal L}_\sigma|_{{\mathcal M}'_\sigma} \cong {\mathcal L}'_\sigma$.  By Hartogs Theorem, it
follows that there is a unique lift of the action of $\Gamma$ to this family
of line bundle $\{{\mathcal L}_\sigma\}_{\sigma\in T}$ over the family of moduli spaces $\{{\mathcal M}_\sigma\}_{\sigma\in T}$ parametrised by Teichm\"{u}ller space $T$.

Geometric quantization states that the Hilbert space which arises when one
wants to quantize ${\mathcal M}$ at level $k$ ($k\in {\Bbb N}$) is
\[\hat{Z}_\sigma = H^0({\mathcal M}_\sigma, {\mathcal L}_\sigma^k).\]
We get this way an association of a vector space $\hat{Z}_\sigma$ to each point
$\sigma\in T$.
The following Theorem constitutes the main result
in the geometric approach to the Witten-Reshetikhin-Turaev quantum representations. 
\begin{theorem}[S.Axelrod, S. Della Pietra \& Witten; N. Hitchin; G. Faltings]\label{HC}
\mbox{ }
\newline
The family of vector spaces $\hat{Z}_\sigma$, $\sigma\in T$ forms a holomorphic vector bundle $\hat{Z}$
over $T$. There is a natural projectively flat mapping
class group invariant connection in the holomorphic bundle $\hat{Z}$ over $T$.
\end{theorem} 

Using techniques
from conformal field theory and the theory of loop groups G. Segal; Tsuchiya,
Ueno \& Yamada (see also \cite{Ts}); A.A. Beilinson \& D. Kazhdan; A.A. Beilison
 \& V. V. Schechtman;
 proved this Theorem in a slightly different formulation, which by
the work of Kumar, Narashimhan \& Ramanathan \cite{KNR} and Beauville \& Laszlo
\cite{BL} is now known to be equivalent to the above stated Theorem. As stated
the Theorem was first proved in the $G=SU(n)$-case by S. Axelrod, S. Della
Pietra \& Witten \cite{ADW} and then by N. Hitchin \cite{Hitchin2} using more 
algebraic
geometric techniques. G. Faltings \cite{F} has generalized N. Hitchin's approach
to the general $G$-case. For a pure differential geometric approach to the Hitchin connection see also \cite{A3}, \cite{AG} and \cite{AGL}.

The mapping class group equivariant cohomology class of the curvature of this connection is known to be non-trivial (see e.g. \cite{Hitchin2} and \cite{TUY}). Hence it is not possible 
to ``correct'' this connection in a mapping class group invariant way so as to obtain a flat connection. 

In order to remedy this non-flatness, we need to tensor
the bundle with a certain fractional power of the determinant bundle over
Teichm\"{u}ller space, which only a central extension of the mapping class group
acts on. It is this change which requires us to introduce some
extra structure called a rigging on the surface $\Sigma$. We refer the reader to
the paper \cite{Walker}, \cite{Segal}, \cite{T} and \cite{AU2} for a thorough explanation of this fact.

Let us now describe the category of closed oriented rigged surfaces. The following
description of this category is identical to the category of extended surfaces
given in K. Walkers paper  \cite{Walker}, except there will be no
piecewise-linear structures here and all surfaces considered will be closed (see also the accounts by V. Turaev in \cite{Tur} and by P. Gilmer and G. Masbaum in \cite{GMas}).
\begin{definition} A rigged surface $(\Sigma,L)$ is a pair consisting of a
closed oriented surface and a Lagrangian subspace $L$ of $H_1(\Sigma,{\Bbb R})$.
\end{definition}
\begin{definition} A rigged morphism from a rigged surface $(\Sigma_1,L_1)$ to a
rigged surface $(\Sigma_2,L_2)$ is a pair $(f,n)$, where $f : \Sigma_1
\rightarrow \Sigma_2$ is an isotopy class of orientation preserving
diffeomorphisms and $n$ is an integer. 
\end{definition}

 We shall now define a
composition of rigged morphisms. Suppose we have two rigged morphisms $(f_1,n_1)
: (\Sigma_1,L_1) \rightarrow (\Sigma_2,L_2)$ and $(f_2,n_2) : (\Sigma_2,L_2)
\rightarrow (\Sigma_3,L_3)$. We then define \[(f_2,n_2)(f_1,n_1) = (f_2f_1, n_2
+ n_1 + \sigma((f_2f_1)_*L_1,(f_2)_*L_2,L_3)),\] where $\sigma(\cdot,\cdot,\cdot)$ is
Wall's signature cocycle for triples of Lagrangian subspaces. (See \cite{Wall}
and also \cite{Walker} Section 18.)

Define the rigged mapping class group $\Gamma_r$ to be the rigged automorphisms
of a rigged surface $(\Sigma,L)$. The rigged mapping class group is a central
extension of the mapping class group of $\Sigma$. The corresponding extension
2-cocycle is the Shale-Weil cocycle. (See Section 17 in \cite{Walker}.) In the
next Section we shall see how the canonical framing choice on 3-manifolds
specified by Atiyah associates to each element of the mapping class group a
preferred element of the rigged mapping class group.  

Recall the construction of the determinant line bundle ${\mathcal L}_D$ over
Teichm\"{u}ller space $T$. For a point $\sigma\in T$ we define \[{\mathcal
L}_{D, \sigma} = \Det H^1(\Sigma_\sigma, {\mathcal O}).\]
There is a natural lift of the action of $\Gamma$ on $T$ to an action on ${\mathcal L}_D$.
As it is described in \cite{Walker} pp. 118-120, the rigged mapping class group
$\Gamma_r$ naturally acts on any real power of the determinant line bundle
${\mathcal L}_D$. We refer the reader to this paper for the explicit
construction.

Having described the category of rigged surfaces let us move to the description
of the functor $Z^{k}_G$. As mentioned in the introduction, the functor $Z^{k}_G$ is a
functor from the category of rigged surfaces to the category of finite
dimensional complex vector spaces. We refer to \cite{Walker}, \cite{Segal} and \cite{AU2}
for the description of the axioms this functor has to satisfy. 
We consider the vector bundle 
$$Z = \hat{Z} \otimes  {\mathcal
L}_{D}^{-\frac{1}{2}\zeta}$$ 
over Teichm\"{u}ller space. Here $\zeta$ is the so
called central charge of the theory given by 
\begin{equation}
\zeta = \frac{|G|k}{k+h},\label{cc}
\end{equation} 
where
$h$ is the dual Coxeter number of $G$, i.e the quadratic Casimir of the adjoint
representation of $G$. (See p.20 in \cite{KNR} for a list of the value of $h$
for all the simple algebras.) The volume $|G|$ is calculated with respect to the
chosen bi-invariant inner product on $G$. 

If $\hat{Z}$ is replaced by the sheaf of vacua in the above Definition of $Z$, one of the main Theorems of \cite{AU2} states that the resulting bundle
has a natural rigged mapping class group invariant connection, which is the tensor product of the TUY-connection in the sheaf of vacua and 
a certain connection in ${\mathcal L}_{D}^{-\frac{1}{2}\zeta}$, which can be constructed from abelian conformal field theory. By the Theorems
of Kumar, Narashimhan \& Ramanathan \cite{KNR} and Beauville \& Laszlo
\cite{BL} and Laszlo \cite{La1}, we know that $\hat{Z}$ with its natural connection from Theorem \ref{HC} above is isomorphic to the sheaf of vacua with its TUY-connection. Hence we conclude that $Z$ has a natural rigged mapping class group invariant flat connection over $T$.
Since $T$ is topologically an open cell of dimension $6g-6$, we can
simply make the following Definition.

\begin{definition} The functor $Z^{(k)}_G$ from the category of genus $g$ rigged
surfaces to the category of finite dimensional complex vector spaces is defined
by associating to each rigged surface $(\Sigma,L)$ the vector space
$Z^{(k)}_G(\Sigma)$ consisting of the global covariant constant sections of $Z$ over
$T$. \end{definition}

In \cite{AU2} we proved using the sheaf of vacua construction model and a number of results of \cite{TUY}, that this gives a modular functor.

From this Definition it is clear that we get formula (\ref{F34}) stated in the introduction.

\section{Atiayh 2-framings of mapping tori and the rigged mapping class group.} 
\label{DEF} 

Following K. Walker we shall now describe how the Atiyah 2-framings of mapping tori gives a set-theoretical section from $\Gamma$ to $\Gamma_r$ (see Chapter 1, ``Extended 2- and 3- manifolds''
p.7-13 and Chapter 17 ``Central Extensions of mapping class groups'' and 18 ``Non-additivity of the signature'' in \cite{Walker}.).
Let $f\in \Diff_+(\Sigma)$ be arbitrary. We need to
 choose a lift $\tilde{f}$ of $f$ in the extended mapping class group of
$\Sigma$. 

Using 2-framings of 3-manifolds
Atiyah defined in $\S$3 of \cite{Atiyah2} a central extension $\widehat{\Gamma}$ of
the mapping class group $\Gamma$ of $\Sigma$. We shall now give an explicit
isomorphism $\Phi$ from $\widehat{\Gamma}$ to $\Gamma_r$ which makes the
following diagram commutative 
\[\begin{array}{ccccccccc} 0 & \rightarrow {\Bbb
Z} & \rightarrow \widehat{\Gamma} & \rightarrow & \Gamma & \rightarrow & 0 \\ &
\mbox{ } \| &  \mbox{ } \mbox{ } \downarrow {\Phi} & & \mbox{ }\|& & \\ 0 & \rightarrow {\Bbb Z} &
\rightarrow \Gamma_r & \rightarrow & \Gamma & \rightarrow & 0 \end{array}\]

Choose a Lagrangian subspace $L\subset H^1(\Sigma,{\Bbb R})$ such that there
exists $M^\dagger$, a compact 3-manifold with the properties that $\partial
M^\dagger = \Sigma$ and $\im(H^1(M^\dagger) \rightarrow H^1(\Sigma)) = L$. Let $M^\dagger \cup_{\id} \Sigma\times I \cup_f
(-M^\dagger)$ be the oriented closed 3-dimensional manifold one obtains by
performing the indicated glueing. (The notation $-M^\dagger$ means $M^\dagger$
with the opposite orientation.) Let $V$ be a compact 4-dimensional manifold such
that $\partial V = M^\dagger \cup_{\id} \Sigma\times I \cup_f (-M^\dagger)$ and
$\sigma(V) = 0$. Here we (also) use $\sigma$ to denote the signature of a
4-manifold. Let $W= M^\dagger\times I$. We note that $\partial W =
M^\dagger \cup_{\id} \Sigma\times I \cup_{\id} (-M^\dagger)$. Identifying the
copy of $M^\dagger$ (resp. $-M^\dagger$) in $\partial W$ with the copy of
$-M^\dagger$ (resp. $M^\dagger$) in $\partial V$, we obtain a compact
4-manifold, denoted $V\cup W$, whose boundary is $\partial(V\cup W) \cong
\Sigma_f$. 

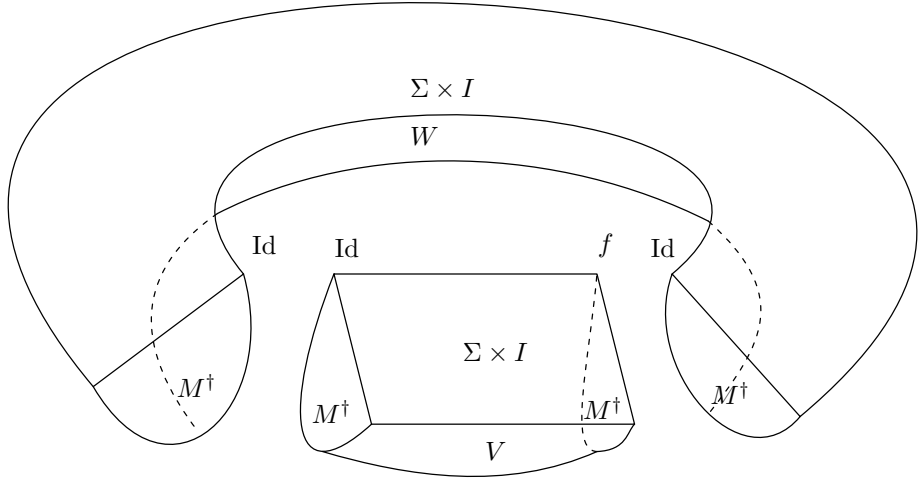
\begin{center}
\begin{texdraw}
\drawdim cm \linewd 0.02
\arrowheadtype t:V 
\move (-2 1) \lvec(1.5 1) \lvec (2 -1) \lvec (-1.5 -1) \lvec (-2 1)
\clvec (-2.7 -.8)(-2.6 -2)(-1.5 -1)
\move (-2.15 -1.365) 
\clvec (-1 -1.7)(.3 -1.9)(1.5 -1.365)
\move (2 -1) \clvec (1.9 -1.1)(1.9 -1.35)(1.5 -1.365)
\lpatt (.067 .1) \clvec (1.1 -1.35)(1.4 -.15)(1.5 1)
\htext(-.3 -.2){$\Sigma\times I$}
\htext(-2 1.2){Id}
\htext(1.5 1.2){$f$}
\htext(0 -1.5){$V$}
\htext(-2.3 -1){$M^\dagger$}
\htext(1.3 -.95){$M^\dagger$}
\lpatt ()
\move(2.5 1) \lvec(4.2 -.9) \clvec (11.85 5.52)(-11.63 7.16)(-5.2 -.5)
             \lvec(-3.2 1) \clvec(-5.77 4.06)(5.56 3.57)(2.5 1)

\clvec(2,05 -.25)(3.43 -1.82)(4.2 -.9) 
\move(-3.2 1) \clvec(-2.7 -.8)(-4.2 -2.2)(-5.2 -.5)

\move(-3.55 1.8) \clvec(-1.6 2.8)(1 2.7)(2.98 1.7)
\lpatt (.067 .1) \clvec(4.1 .7)(3.7 .05)(2.98 -.86)
\move(-3.55 1.8) \clvec(-4.7 1)(-4.65 0)(-3.8 -1.1)
\lpatt()
\htext(-4.1 -.65){$M^\dagger$}
\htext(3 -.75){$M^\dagger$}
\htext(-1 3.3){$\Sigma\times I$}
\htext(-1 2.7){$W$}
\htext(-3.1 1.25){Id}
\htext(2.2 1.2){Id}

\move(0 5)
\end{texdraw}
\end{center}

\bigskip
\bigskip

According to Wall's signature Theorem in \cite{Wall}, we have that
\[\sigma(V\cup W) = \sigma_L(f),\] where we have used the notation \[\sigma_L(f)
= \sigma(L\oplus L,\graph(-f),\graph(-id)).\] Suppose now $\alpha$ is a
2-framing of $\Sigma_f$, i.e. $(\Sigma_f,\alpha)\in \widehat{\Gamma}$. We then
define \[{\Phi}(\Sigma_f,\alpha) = (f, \frac{1}{6}p_1(2T_{W\cup V},\alpha)).\]
which gives a group homomorphism from $\widehat{\Gamma}$ to $\Gamma_r$
that obviously makes the above diagram commutative.

From the definition of ${\Phi}$ we see that the mapping torus $\Sigma_f$ with its
Atiyah 2-framing $S(f)\in \widehat{\Gamma}$ is by ${\Phi} : \widehat{\Gamma}
\rightarrow \Gamma_r$ mapped to $\tilde{f} := (f,\sigma_L(f))\in \Gamma_r.$

\section{Seifert invariants of finite order mapping tori.} 
\label{SI} 

The  three manifold
type of the mapping torus $\Sigma_f$ is determined by the Nielsen-Thurston type
of the diffeomorphism $f$.  We are here only interested in the diffeomorphisms
which preserve a point in Teichm\"{u}ller space. This is of course equivalent to
requiring that the diffeomorphism $f$ is of finite order.
Since $f$ is finite order, $\Sigma_f$ is a Seifert fibred 3-manifold.
Let us here review how the Seifert invariants of $\Sigma_f$ are related to $f$.
We refer to \cite{O} for the definition of the Seifert invariants of a Seifert
fibred manifold. (See pp. 11-16, Theorem 3 p. 90 and Theorem 6 p. 97 in \cite{O}.)
Let $\sigma\in T$ be a fixed point for $f$. We then have that $f$ is an
automorphism of the Riemann surface $\Sigma_\sigma$. Let ${\tilde\Sigma} =
\Sigma/f$ and let $\pi : \Sigma \rightarrow {\tilde\Sigma}$ be the
projection. Then $\Sigma$ is an $m$-fold branched cover of the Riemann Surface
${\tilde\Sigma}$. The generic point ${\tilde p}\in {\tilde \Sigma}$ has
the property that $\# \pi^{-1}({\tilde p}) = m$. Let
$\{{\tilde p}_1,\ldots,{\tilde p}_n\}$ be the finite set of special points
on ${\tilde \Sigma}$ for which $\#\pi^{-1}({\tilde p}_i) = m_i < m$. Assume
that $m_i = 1$ for $i = 1, \ldots, l$ say (in the case $f$ has no fixed points we set $l=0$). Let $p_i\in \pi^{-1}({\tilde p}_i).$
Then we have that 
\[\pi^{-1}({\tilde p}_i) = \{ p_i,f(p_i),\ldots,
f^{m_i-1}(p_i)\},\]
 and $f^{m_i}(p_i) = p_i$. From this we see that $m_i|m.$ Let
$l_i$ be defined by $m_il_i = m.$

Since $f^{m_i}(p_i) =p_i$ there is a local coordinate $\xi_i$ around $p_i$ such
that 
\[\xi_i\circ f^{m_i}  = \expo (2\pi i \frac{n_i}{l_i}) \xi_i,\] 
for some $0<n_i<l_i$
satisfying $(n_i,l_i) = 1$. Let $0<k_i<l_i$ be given by $k_in_i = 1 \mod l_i$.
This uniquely determines $k_i$. The Seifert invariants in the notation of
\cite{O} of the mapping torus $\Sigma_f$ is (see p.90 in \cite{O}) 
\[
(b,g,(\alpha_1,\beta_1),\ldots, (\alpha_n,\beta_n))(\Sigma_f) = (-\sum_{i=1}^n
\frac{l_i}{k_i},g({\tilde \Sigma}), (l_1,k_1),\ldots,(l_n,k_n)),\] 
and
$\Sigma_f$ is of type $o_1$, i.e. $\epsilon = o_1$ (see p. 88 in \cite{O}).

\begin{definition} The orbifold Euler number of the Orbifold $S^1$ bundle which
constitutes a Seifert fibred manifold $M$ with the Seifert invariants
$(b,g,(\alpha_1,\beta_1),\ldots, (\alpha_n,\beta_n))$ is called the Seifert
fibred Euler number and it is given by \[e(M) = - (b + \sum_{i=1}^n
\frac{\beta_i}{\alpha_i}).\]

\end{definition}

\noindent Notice that the orbifold Euler number of $\Sigma_f$ is \[e(\Sigma_f) =
0.\]

Let us now analyse which Seifert fibred manifolds are diffeomorphic to the
mapping torus of some finite order diffeomorphism of a closed oriented surface.
This can be done by simply combining a few of the Theorems in \cite{O} and
\cite{Scott}. 

\begin{theorem}[Orlik] A Seifert fibred manifold $M$ of type $o_1$ and with Seifert
invariants $(b,g,(\alpha_1,\beta_1),\ldots, (\alpha_n,\beta_n))$ is
diffeomorphic to the mapping torus of a finite order diffeomorphism of a closed
oriented surface iff \[e(M) = 0.\] \end{theorem}

\proof The Theorem for Large Seifert fibred manifolds follow
directly from Theorem 5.4(ii) p. 167 in \cite{Scott} and Corollary 5 p. 122 in
\cite{O}. One then looks through the list of small Seifert fibred manifolds of
type $o_1$ which fibres over the circle (p.124 \cite{O}) and checks that the
Theorem holds for these 4 cases too. \begin{flushright} $\Box$ \end{flushright}

A simple calculation shows that there are many Seifert fibred
manifolds with $e=0.$ Let $n$ pairs of integers $(\alpha_1,\beta_1),\ldots,
(\alpha_n,\beta_n)$, such that $(\alpha_i,\beta_i) = 1$ and $0<\beta_i<\alpha_i$
be given. Then one can choose $(\alpha_{n+1},\beta_{n+1})$ such that there is 
a Seifert fibred manifold with
$e=0$ and the Seifert invariants $(b,g,(\alpha_1,\beta_1),\ldots,
(\alpha_{n+1},\beta_{n+1}))$ for any non-negative integer $g$.

\section{The framing correction.} 
\label{FC}

We shall in this Section compute the framing correction $ \Det(f)^{-\frac{1}{2}\zeta}$, hence we choose a Lagrangian
subspace $L\subset H^1(\Sigma,{\Bbb R})$, calculate 
$$\sigma_L(f) = \sigma(L\oplus L,
\graph(-f),\graph(-id))$$ 
and then calculate $f$'s action on ${\mathcal
L}_D^{-\frac{1}{2}\zeta}$ using $L$. Based on this we calculate the framing
correction $\Det(f)^{-\frac{1}{2}\zeta}$.
Let us first choose $L$ and calculate $\sigma_L(f)$. Since $f$ preserves the complex
structure $\sigma$, we just need to compute $f$'s action on $H^0(\Sigma_\sigma,\Omega^1) \cong
H^1(\Sigma,{\Bbb R}).$
Since $f^m=1$, we know that there exists a complex basis $(e_1, \ldots ,e_g)$ for
$H^0(\Sigma_\sigma,\Omega^1)$ such that 
\[f e_i = \omega_i e_i,\] 
where each $\omega_i$ is an m'th-root of unity.
Choose \[L=\spane_{\Bbb R}\{e_1,\ldots,e_g\}\] and notice that $L$ is a Lagrangian
subspace.

We now consider the symplectic vector space $V = H^1(\Sigma,{\Bbb R}) \oplus
(-H^1(\Sigma,{\Bbb R})).$ (The minus sign indicates that we are considering
the sum of the symplectic structure on the first copy and minus the symplectic
structure on the second copy.) Let \[L_1 = L\oplus L \subset V.\] 
Define 
\[L_2 = \{(x,-fx) | x \in H^1(\Sigma,{\Bbb R})\}\] 
and 
\[L_3 = \{(x,-x) | x \in H^1(\Sigma,{\Bbb R})\}.\] 
If we have that $\omega_i = \omega^R_i + i \omega^I_i$ and
we use the real basis $(e_1, i e_1, \ldots, e_g, ie_g)$, then we get that 
\begin{eqnarray*}L_1
& = & \spane_{\Bbb R} \{ (e_1,0), \ldots ,(e_g,0), (0,e_1), \ldots ,(0,e_g)\},\\
L_2&  = &  \spane_{\Bbb R} \{ (e_1,-\omega^R_1 e_1 - \omega^I_1 i e_1), \ldots
,(e_g,-\omega^R_g e_g - \omega^I_g i e_g), \\
& & \mbox{ } \mbox{ }\mbox{ }\mbox{ }\mbox{ }\mbox{ }\mbox{ }\mbox{ }\mbox{ }\mbox{ }(ie_1,\omega^I_1 e_1 - \omega^R_1 ie_1), \ldots
,(ie_g,\omega^I_g e_g - \omega^R_g ie_g)\}\\
L_3 & = & \spane_{\Bbb R} \{ (e_1,-e_1),
\ldots ,(e_g,-e_g), (ie_1,-ie_1), \ldots ,(ie_g,-ie_g)\}.
\end{eqnarray*}
Consider the
following symplectic subspace of $V$ \[W_i = \spane_{\Bbb R} \{ (e_i,0),
(ie_i,0) ,(0,e_i), (0,ie_i)\},\] and set $L_j^i = L_j\cap W_i.$
We wish to decompose $W_i$ into a direct sum of 2 symplectic subspaces
 $W_i = W'_i \oplus W''_i,$ such that $W'_i\cap L^i_j$ and $W''_i\cap
L^i_j$ are Lagrangian. The case $\omega^R_i \neq 1$ will be considered first. Let
\[W'_i = \spane_{\Bbb R} \{ (e_i,e_i), (ie_i, -ie_i)\}\] and \[W''_i =
\spane_{\Bbb R} \{ (e_i,-e_i), (ie_i,\omega^I_i e_i -\omega^R_i ie_i)\}.\]
We notice that $W_i = W'_i \oplus W''_i,$ and both $W'_i$ and $W''_i$ are
symplectic subspaces. We get that 
\[L_1\cap W'_i = \spane_{\Bbb R} \{(e_i,e_i)\},\] 
\[L_2\cap W'_i = \spane_{\Bbb R} \{(\omega^R_i-1)(e_i,e_i)-\omega^I_i(ie_i,-ie_i)\},\] 
\[L_3\cap W'_i = \spane_{\Bbb R} \{(ie_i,-ie_i)\}\] 
and 
\[L_1\cap W''_i = L_3\cap W''_i = \spane_{\Bbb R} \{(e_i,-e_i)\},\] 
\[L_2\cap W''_i = \spane_{\Bbb R} \{(ie_i,\omega^I_ie_i-\omega^R_i ie_i)\}.\] 
The case $\omega^R_i = 1$ is trivial since it gives $L_2^i =
L_3^i.$
We can now calculate 
\[\sigma(L_1,L_2,L_3)= \sum_{i=1}^g \sigma_{W_i}(L_1\cap
W_i, L_2\cap W_i, L_3\cap W_i).\] 
Because of coinciding subspaces $\sigma_{W''_i}(L_1\cap W''_i, L_2\cap W''_i,
L_3\cap W''_i) = 0$, so we just need
to calculate $\sigma_{W'_i}(L_1^i\cap W'_i, L_2^i\cap W'_i, L_3^i\cap W'_i).$ We note that
$((e_i,e_i),(ie_i,-ie_i))$ is a positive oriented basis for $W'_i$ so
\begin{eqnarray*} \sigma_{W'_i}(L_1^i, L_2^i, L_3^i) & = & - \sign
\frac{\omega^R_i-1}{\omega^I_i}  \\ & = &
\sign \, \im \omega_i. \end{eqnarray*} Hence we get the formula
\begin{equation}\sigma(L\oplus L, \graph(-f),\graph(-id)) = \sum_{i=1}^g \sign \,
\im \omega_i.\label{AF} \end{equation}
Let us now calculate $\Det(f)^{-\frac{1}{2}\zeta}$. From the expression of $f$'s
action on $H^0(\Sigma_\sigma,\Omega^1)$ we see that 
\[\Det(f) = \prod_{i=1}^g
\omega_i\] and \[f(L) = \spane_{\Bbb R}\{\omega_i e_i\}.\] 
Suppose that
$\omega_i= \pm 1$, for $i= 1,\ldots,\iota$ and that $\omega_i\neq \pm 1$, for $i=
\iota+1,\ldots,g$. Then let 
\[\begin{array}{lclr} f_i^{\pm}& = & e_i & 1\leq i \leq \iota\\
 f_i^\pm & = & \pm(\frac{\omega^R_i}{\omega^I_i}e_i + i e_i) & \iota+1\leq i \leq g.
\end{array}\] 
The basis $(f^\pm_1,\ldots,f_g^\pm)$ is a basis for $f(L)$
satisfying the condition on p. 119 in Chapter 16 in \cite{Walker}. We then consider
\[P_t^\pm(L,f(L)) = \spane_{\Bbb R}\{(1-t)e_i + t f^\pm_i\}.\] 
Now notice that
$P_t^\pm(L,f(L))$ is homotopic to $\tilde{P}_t^\pm(L,f(L))$ when we let
\[\tilde{P}_t^\pm(L,f(L)) = \spane_{\Bbb R}\{\omega_i^{t,\pm} e_i\},\] 
where for $t\in [0,1]$
\[\omega_i^{t,\pm} = \left\{ \begin{array}{cc} 1 & \omega_i = \pm 1.\\ \left\{
\begin{array}{cc} \exp(it \log\omega_i) & \pm = +\\ \exp(it(-\pi +
\log\omega_i)) & \pm = - \end{array}\right. & \im \omega_i > 0\\ \left\{
\begin{array}{cc} \exp(it \log\omega_i) & \pm = +\\ \exp(it(\pi + \log\omega_i))
& \pm = - \end{array}\right. & \im \omega_i < 0 \end{array}\right. \]
 Here
$\log\omega_i \in (-\pi,\pi).$ We then have that the path
$\Det^2_{(1-t)e+tf^\pm}(P_t^\pm(L,f(L)))$ is homotopic rel. end points to
$\Det^2_{\omega^{t,\pm}e}(\tilde{P}_t^\pm(L,f(L)))$. Clearly
\[\Det^2_{\omega^{t,\pm}e}(\tilde{P}_t^\pm(L,f(L))) = (\prod_{i=1}^g
\omega_i^{t,\pm})^2 (e_1\wedge\ldots\wedge e_g)^2.\]

Let $\tilde{\omega}_i^\pm$ be the unique curve in ${\Bbb C}$ which projects under the exponential map $\expo(i\cdot)$ to
$(\omega_i^{t,\pm})^2$ and starts in $0$. Explicitly, we get for $t\in [0,1]$ that
\[\tilde{\omega}_i^{t,\pm} = \left\{ \begin{array}{cc} 1 & \omega_i = \pm 1.\\
\left\{ \begin{array}{cc} 2t \log\omega_i & \pm = +\\ 2t(-\pi + \log\omega_i) &
\pm = - \end{array}\right. & \im \omega_i > 0 \\ \left\{ \begin{array}{cc} 2t
\log\omega_i & \pm = +\\ 2t(\pi + \log\omega_i) & \pm =- \end{array}\right. &
\im \omega_i < 0 \end{array} \right.  \] 
By the description of the action of
$\Gamma_r$ on $\tilde{{\mathcal L}}_D$, we get that $(f,n)$ acts on
$\tilde{{\mathcal L}}_{D,\sigma}$ by 
\[(f,n)(p) = p +
\left(\sum_{\stackrel{i=1}{\omega_i\neq \pm 1}}^g 2\log \omega_i + \pi (-
\#\{\omega_i|\im\omega_i > 0\} + \#\{\omega_i|\im\omega_i<0\} + n)\right)\]
for $p\in \tilde{{\mathcal L}}_{D,\sigma}$.
Here $+$ refers to the action
of ${\Bbb C}$ on $\tilde{{\mathcal L}}_{D,\sigma}$. From this expression, we see
that the correction term coming from the choice of the Atiyah 2-framing, i.e. 
$$n=\sigma_L(f) = \sum_{i=1}^g \sign \, \im \omega_i,$$ 
cancels with the 
above multiple of $\pi$ and we get the simple
formula 
\[\Det(f)^{-\frac{1}{2}\zeta} = \prod_{\stackrel{i=1}{\omega_i\neq \pm
1}}^g \expo(-i \frac{1}{2}\zeta \log \omega_i),\] 
where $\log\omega_i \in
(-\pi,\pi)$ for $\omega_i\neq \pm 1.$

Let us now apply the Lefschetz fixed point formula to get an expression for the
framing correction in terms of the Seifert invariants of the mapping torus
$\Sigma_f$. Since the cyclic group $\langle f\rangle$ of order $m$ acts on
$H^0(\Sigma_\sigma,\Omega^1)$, there is a unique decomposition 
\[ H^0(\Sigma_\sigma,\Omega^1) =
\bigoplus_{\alpha^m = 1} M_\alpha,\] 
where $M_\alpha$ is the $\alpha$-eigenspace
of $f$. Let $d_\alpha = \dim_{\Bbb C} M_\alpha,$ and denote by $\tr(f^\beta)$ the
trace of $f^\beta$'s action on $H^0(\Sigma_\sigma,\Omega^1)$. Note that
$\Det(f)^{-\frac{1}{2}\zeta}$ only depends on the $d_\alpha$'s since
\[\Det(f)^{-\frac{1}{2}\zeta} = \prod_{\stackrel{\alpha^m=1}{\alpha\neq \pm 1}}
\alpha^{-\frac{1}{2}\zeta d_\alpha}.\] 
(This fractional power should be
interpreted with respect to the above branch of the logarithm.) One easily
verifies that \[d_\alpha = \frac{1}{m} \sum_{\beta = 0}^{m-1} \tr(f^\beta)
\alpha^{-\beta}.\] 
The Lefschetz-Riemann-Roch Theorem due to Atiyah \& Bott now
gives 
\[1-\tr(f^\beta) = \sum_{j=1}^l (1-\expo(2 \pi i \frac{n_j \beta}{m}))^{-1}
\mbox{, } \beta \neq 0 \mod m\] 
and by the very definition of the genus of a
Riemann surface, we have \[\tr(f^0) = g.\] 
From this we then get 
\[ d_\alpha =
\frac{1}{m} \left( - \sum_{j=1}^l \sum_{\beta = 1}^{m-1}
\frac{\alpha^{-\beta}}{1-\expo(2\pi i \frac{n_{j} \beta}{m})} + g + \sum_{\beta
= 1}^{m-1} \alpha^{-\beta}\right) .\] 
But since
\[\sum_{\beta= 0}^{m-1} \alpha^{-\beta}=0\]
for $\alpha\neq  1, \alpha^m=1,$ we see that $\sum_{\beta
= 1}^{m-1} \alpha^{-\beta} = -1$. 
Define for $\alpha^m = 1$, and $(n,m)=1$, $m> 1$
\[\mu_\alpha(n) = - \sum_{\beta =
1}^{m-1} \frac{\alpha^{-\beta}}{1-\expo(2\pi i \frac{n \beta}{m})}.\]
We have the following elementary claim, which computes $\mu_\alpha(n)$
\begin{claim}
Let $\overline{n}\in {\Bbb Z}$ be the unique integer determined by 
\[\alpha = \expo(2\pi i\frac{n\overline{n}}{m}), \mbox{ } 0\leq\overline{n}<m.\]
Then 
\[\mu_\alpha(n) = \overline{n} - \frac{m-1}{2}.\]
\end{claim}
One can prove this claim by an elementary calculus proof, by multiplying the 
last term in the denominator of the expression for $\mu_\alpha(n)$ by a free
variable $x$ say and then examining the resulting function of $x$ near $x=1$,
using the geometric series expansion of the resulting function.
We shall introduce the following notation for later use
\[\mu_m^i(n) = \mu_{e^{2\pi i \frac{i}{m}}}(n).\]
From the above calculation we get that
\begin{theorem}\label{fct}
The framing correction term is given by
 \begin{equation}\Det(f)^{-\frac{1}{2}\zeta} =
\prod_{\stackrel{\alpha^m=1}{\alpha\neq \pm 1}} \alpha^{-\frac{1}{2}\zeta
\frac{1}{m}\left( g-1+\sum_{j=1}^l \mu_\alpha(n_j) \right) } .\label{FCT}
\end{equation}
\end{theorem}
This formula expresses the framing correction completely in terms of the
Seifert invariants of $\Sigma_f$ and $m$.
\begin{remark}
When we specialize to the case where $f$ has no fixed points, the framing correction term given by (\ref{FCT}) reduces to
\[\Det(f)^{-\frac{1}{2} \zeta} = \prod_{\stackrel{\alpha^m=1}{\alpha\neq \pm 1}}
\alpha^{-\frac{1}{2}\zeta (g(\tilde{\Sigma})-1)} = 1.\]
\end{remark}

\section{The fixed point set.} 
\label{FPS} 
Recall that ${\mathcal M}$ denotes the moduli space of flat $G$-bundles over $\Sigma$.
We want to describe the fixed point
set $|{\mathcal M}|\subset {\mathcal M}$ for the finite order automorphism $f$ of $\Sigma$. We shall
do this in the general case, where we are considering a semi-simple, simply
connected Lie group $G$. Let $P$ be a principal $G$-bundle over $\Sigma$ and ${\mathcal G}$ the gauge group of $P$. Since
the isomorphism classes of principal $G$-bundles are classified by homotopy
classes of maps from $\Sigma$ to $BG$, which in this particular case are in 1-1 
correspondence
with $H^2(\Sigma,\pi_1(G)) = \{0\},$ we see that $P$ is topological trivial.

Suppose $A$ is a connection in $P$ which represents a fixed point for $f$ in ${\mathcal M}.$ Since $f^*P$ is isomorphic to $P$ and $[f^*A]=[A]$ in ${\mathcal M}$, we know there is an isomorphism $\psi$ from $P$ to $f^*P$ such that $\psi f^*A=A$.
We notice that when we compose $\psi$ with the natural bundle map from $f^*P$ to
$P$ covering $f$, we get a lift ${\phi}: P \rightarrow P$ covering $f$. Since $f^m = 1$, we see that ${\phi}^m$ covers the 
identity and ${\phi}^m A=A$, so ${\phi}^m\in Z_A,$ the stabilizer of $A$ in ${\mathcal G}$.
From this we see that for any connection $A$ representing a fixed point in ${\mathcal M}$, we get a lift ${\phi}$ such that  ${\phi} A=A$, but we only get that 
${\phi}^m\in Z_A$ and not necessarily that ${\phi}^m=1.$

We shall first concentrate on describing the part of the fixed point set, which is
contained in the subset of irreducible connections $|{\mathcal M}'|\subset
{\mathcal M}'.$ Recall that ${\mathcal M}'$ is by definition the subset consisting of equivalence classes of flat connections which have stabilizer equal to the center $Z(G)$.

Now consider the set $L$ of lifts ${\phi}$ of $f$ such that $ {\phi}^m \in
Z=Z(G).$  For each $\phi\in L$ let ${\mathcal A}'_{{\phi}}$ be the space of flat irreducible
connections in $P$ invariant under ${\phi}$ and let
$${\mathcal G}_{{\phi} }= \{ g \in {\mathcal G} \mid [g,{\phi}]  \in Z\}.$$ 
We observe that ${\mathcal G}_{{\phi} }$ acts on ${\mathcal A}'_{{\phi}}$. Moreover, if $g\in {\mathcal G}$ has the property that
$$g({\mathcal A}'_{{\phi}}) \cap {\mathcal A}'_{{\phi}} \neq \emptyset,$$ then $g\in {\mathcal G}_{{\phi} }$.

We have an action of $Z\times {\mathcal G}$ on $L$ given by 
$$(z,g)({\phi}) = z g {\phi} g^{-1}.$$
Let us denote the set of equivalence classes of
lifts of $f$ by $\Delta = L/(Z\times {\mathcal G})$.
We note that if $(z,g)({\phi}_1) = {\phi}_2$, then $g$ induces a bijection between ${\mathcal A}'_{\phi_1}$ and
${\mathcal A}'_{\phi_2}$. 

Suppose there exist a $g\in {\mathcal G}$ such that $g({\mathcal A}'_{\phi_1}) \cap {\mathcal A}'_{\phi_2} \neq \emptyset.$ 
Then let $A_i\in {\mathcal A}'_{\phi_i}$, $i=1,2$, such that $gA_1 = A_2$. Now we compute that
$$g \phi_1 g^ {-1} A_2 = \phi_2 A_2,$$
so there exists $z\in Z$ such that 
$$ z g \phi_1 g^ {-1} = \phi_2.$$

From the above discussion we conclude that
\begin{lemma}

We have the following description of the fixed point set inside the moduli space of irreducible connections
$$ |{\mathcal M}'| = \coprod_{[{\phi}]\in \Delta} {\mathcal A}'_{{\phi}}/{\mathcal G}_{{\phi} }.$$

\end{lemma}

Let us define for $\delta\in \Delta$
$$ |{\mathcal M}'|_\delta = {\mathcal A}'_{{\phi}}/{\mathcal G}_{{\phi} } \subset |{\mathcal M}'| ,$$
for any ${\phi} \in L$, which represents $\delta$.

First we determine the finite set $\Delta$ explicitly. Let $C$ be the set of conjugacy classes of $G$. We have a natural map
$$ H : \Delta \rightarrow  \{(z,c_1, \ldots c_n) \in Z \times C^n\mid c_i^{l_i} = z\}/Z$$
induced by mapping ${\phi}$ to $({\phi}^m, {\phi}^{m_1}(p_1), \ldots, {\phi}^{m_n}(p_n)).$ The action of $Z$ is as follows
$$z'(z,c_1, \ldots, c_n) = ((z')^m z, (z')^{m_1}c_1, \ldots, (z')^{m_n}c_n).$$

\begin{prop} \label{H}
The map $H$ is a bijection.
\end{prop}

Before we get into the proof of this Proposition, let us introduce the following notation. Choose small embedded closed discs $D_i$ around each of the special points $p_i$, $i=1, \ldots, n$, and an embedded closed disc $D_0$, such that
$f^j(D_i)$, $j=0, \ldots, m_i-1,$ $i=0,\ldots n$ are all disjoint (we set $m_0 = m$). Let $\Sigma'$ be the compliment of the interior of all these discs. Let ${\tilde \Sigma}' = \Sigma'/ \langle f\rangle$. The component of the boundary of ${\tilde \Sigma}'$ which equals $\pi(\partial D_i)$ we denote $\partial_i({\tilde \Sigma}')$. We further let $\partial_0({\tilde \Sigma}') = \pi(\partial D_0)$. Let us denote by $\bar  G$ the factor group $G/Z$ and by $\bar P$ the $Z$ quotient of $P$, i.e $\bar P$ is the trivial $\bar  G$-bundle over $\Sigma$. A lift $\phi$ of $f$ to $P$ induces a lift $\bar\phi$ to $\bar P$. We further use the following notation for the restrictions
$$P' = P|_{\Sigma'}, \mbox{ } \bar P' = \bar P|_{\Sigma'},$$
and for the lifts
$$ \phi' = \phi|_{\Sigma'}, \mbox{ } \bar \phi' = \bar \phi|_{\Sigma'}.$$
We define $\bar P'_{\bar \phi'} = \bar P'/\langle \bar \phi'\rangle$, which is a principal $\bar G$-bundle over ${\tilde \Sigma}' $. We note that since ${\tilde \Sigma}' $ is not closed, $\bar P'_{\bar \phi'} $ is trivializable as a principal $\bar G$-bundle over ${\tilde \Sigma}' $.
\proof
Let us first prove injectivity. Assume that we have lifts $\phi_i \in L$ such that 
$$H(\phi_1) = H(\phi_2).$$
Let us denote the $z'\in Z$ the element which translates the conjugacy classes of $\phi_2$ to those of $\phi_1$. 
There are unique maps $\varphi_i : \Sigma \ra G$, such that $\phi_i(x,g) = (f(x), \varphi_i(x)g)$, for all $x\in \Sigma$ and $g\in G$. Let $\bar\varphi_i : \Sigma \ra \bar G$ be the composite of $\varphi_i$ and the project from $G$ to $\bar G$.
Choose trivializations 
$$T_i  : \bar P'_{{\bar \phi'}_i}|_{{\tilde \Sigma}'}\ra {\tilde \Sigma}' \times \bar G.$$
Fix isomorphisms
$$ \bar \Psi_i : \bar P'\ra \pi^*\bar P'_{{\bar \phi}'_i}.$$
This way we get a commutative diagram
\[ \begin{CD} \bar P' @>{\bar\Psi_i}>> \pi^*\bar P'_{{\bar \phi}'_i} @>{\pi^*T_i}>> \bar P'  \\
 @V{\bar \phi'_i}VV @V{f\times \Id}VV @VV{f\times \Id}V \\ 
 P' @>{\bar\Psi_i}>> \pi^*\bar P'_{{\bar \phi}'_i} @>{\pi^*T_i}>> \bar P'  
 \end{CD}\]
Let 
$$\bar \Gamma_i = \pi^*T_i \circ \bar \Psi_i$$ and let $\bar \gamma_i : {\Sigma}' \ra \bar G$ be such that
$\bar \Gamma_i(x, \bar g) = (x, \bar \gamma_i(x) \bar g).$ We get the relation 
$$\bar\varphi_i|_{\Sigma'} = (\bar \gamma_i\circ f)^{-1} \bar \gamma_i.$$
Since $(\bar \varphi_i|_{\Sigma'})_* : \pi_1(\Sigma') \ra Z$ is trivial, we get that
$$(\bar \gamma_i)_* \circ f_* = (\bar \gamma_i)_* : \pi_1(\Sigma') \ra Z.$$
Let $\bar g' = \bar \Gamma_2^{-1} \Gamma_1 : \bar P' \ra \bar P' $ and $\bar \gamma' : \Sigma' \ra \bar G$ be given by
$$\bar g'(x,\bar g) = (x, \bar \gamma'(x) \bar g).$$
 Then
\begin{equation}
\begin{CD} \bar P' @>{\bar g'}>>  \bar P'  \\
 @V{\phi'_1}VV  @VV{\phi'_2}V \\ 
 \bar P' @>{\bar g'}>> \bar P'
 \end{CD}\label{g'com}
 \end{equation}
We of course also have that
$$(\bar \gamma')_* \circ f_* = (\bar \gamma')_* : \pi_1(\Sigma') \ra Z.$$
Since $\Sigma' \ra {\tilde \Sigma}'$ is an unramified covering, we have a short exact sequence
$$ 1 \ra \pi_1(\Sigma') \ra \pi_1(\tilde \Sigma') \ra {\Bbb Z}_m \ra 1.$$
This gives rise to the following exact sequence in cohomology (see VII.6.4. in [Br])
$$ 0 \ra H^1( {\Bbb Z}_m , Z) \ra H^1(\pi_1(\tilde \Sigma'), Z) \ra H^1(\pi_1(\Sigma'), Z)^{ \pi_1(\tilde \Sigma')} \stackrel{\delta}{\ra} H^2({\Bbb Z}_m , Z).$$
So $(\bar \gamma')_* \in H^1(\pi_1(\Sigma'), Z)^{ \pi_1(\tilde \Sigma')}$ as demonstrated above. The obstruction of extending $(\bar \gamma')_*$ to a homomorphism form
$ \pi_1(\tilde \Sigma')$ to $Z$ is measured by the element $\delta((\bar \gamma')_*) \in H^2({\Bbb Z}_m , Z) \cong Z/Z^m$. To compute this element of $H^2({\Bbb Z}_m , Z)$ we consider a curve $\epsilon$ on $\Sigma'$, which connects a point, say $\epsilon_0$, and its image under $f$, say $\epsilon_1=f(\epsilon_0)$. Choose a lift $\gamma' : \epsilon \ra G$ of $\bar \gamma'\mid_\epsilon : \epsilon \ra \bar G$ and let $g'$ be the corresponding lift of $\bar g' : \bar P'\mid_{\epsilon} \ra \bar P'\mid_{\epsilon}$ to $P'\mid_{\epsilon}$. Then there exists $\tilde z \in Z$ such that 
$$\tilde z \phi_2' g'|_{\epsilon_0} = g' \phi_2'|_{\epsilon_0}.$$
By abuse of notation, we will now denote $\tilde z \phi_2'$ by just $\phi_2'$. This does of course not change the equivalence class of $\phi_2'$. We now extend $g'$ to $P'|_{f^j(\epsilon)}$, $j=1, \ldots, m-1$ by requiring that the diagram
\[ \begin{CD} P'|_\epsilon @>{g'}>>  P'|_\epsilon \\
 @V{(\phi'_1)^j}VV  @VV{(\phi'_2)^j}V \\ 
 P'|_{f^j(\epsilon)} @>{g'}>> P'|_{f^j(\epsilon)}
 \end{CD}\]
is commutative. This defined $\gamma' : \coprod_{j=0}^{m-1} f^j(\epsilon) \ra G$ such that $\gamma'(f^j(\epsilon_0)) = \gamma'(f^{j+1}(\epsilon_0))$, $j=0,\ldots, m-2,$ and there exists $z''\in Z$ such that 
$$z'' \gamma'(\epsilon_0) = \gamma'(f^{m-1}(\epsilon_1)).$$
By the  construction of $\delta$, we have that
$$\delta((\bar \gamma')_*) = z''$$
under the isomorphism $H^2({\mathbb Z}_m,Z) \cong Z/Z^m.$ But 
\[ \begin{CD} P'|_{\epsilon_0} @>{g'(\epsilon_0)}>>  P'|_{\epsilon_0}   \\
 @V{\phi'_1}VV  @VV{\phi'_2}V \\ 
 P'|_{f^{m-1}(\epsilon_1)}  @>{g'(f^{m-1}(\epsilon_1))}>> P'|_{f^{m-1}(\epsilon_1)}
 \end{CD}\]
hence 
$$(z')^m = \frac{(\phi'_1)^m}{(\phi_2')^m} = \frac{\gamma'(f^{m-1}(\epsilon_1))}{\gamma'(f^{m-1}(\epsilon_0))} = z''.$$
Hence $\delta((\bar \gamma')_*) = 1 \in H^2({\mathbb Z}_m,Z).$ Thus we can extend $(\bar \gamma')_*$ to a homomorphism from $\pi_1(\bar \Sigma')$ to $Z$. This mean 
we can find $\bar g'' \in {\mathcal G}(\tilde\Sigma'\times \bar G)$ such that $(\bar g'')_* \circ \pi_* = (\bar \gamma')_*$. Now redefine $T_2$ to be $\bar g'' \circ T_2$. The resulting new $\bar g'$ still has the property that the diagram (\ref{g'com}) commutes and its associated $(\bar \gamma')_* : \pi_1(\Sigma') \ra Z$ is trivial, so we can lift $\bar g'$ to an element $g'\in {\mathcal G}(P')$, such that 
\begin{equation}
\begin{CD} P' @>{g'}>>  P'  \\
 @V{\phi'_1}VV  @VV{z_0'\phi'_2}V \\ 
 P' @>{g'}>> P'
 \end{CD}\label{g'comz}
\end{equation}
commutes for some element $z_0' \in Z$.

From this we conclude that $(z_0'/z')^m=1.$ We observe now that by changing the choice of $\bar g''$, such that $(\bar g'')_*$ takes a generator of $\pi(\tilde \Sigma')/\pi_*(\pi_1(\Sigma'))$ to $z_0'/z'$, then we get a new $g'$ with the property that the $z_0'$ that appears in (\ref{g'comz}) is equal to $z'$. We can therefore from now on assume this to be the case.

We will now extend $g'$ to an element $g\in {\mathcal G}$, such that (\ref{g'comz}) continuous to be commutative with $P'$ replaced by $P$, $\phi'_i$ replaced with $\phi_i$ and $g'$ by $g$.

First we extend $g'$ over $D_0$ in any way we like. Then we extend $g'$ over $f^j(D_0)$ by the formula
$$ g'(f^j(x)) = z' \varphi_2(f^{j-1}(x)) g'(f^{j-1}(x)) \varphi_1(f^{j-1}(x))^{-1}.$$
From this we get that 
$$g'(f^m(x)) = (z')^m \left( \prod_{i=0}^{m-1} \varphi_2(f^{m-1-i}(x))\right) g'(x) \left( \prod_{i=0}^{m-1}\varphi_1(f^{m-1-i}(x))\right) ^{-1} = g'(x).$$
Hence $g'$ is now also defined and satisfies (\ref{g'comz}) on $f^j(D_0)$, $j=0,\ldots, m-1$.

We observe that $f^{m_i}(D_i) = D_i$, hence we first consider the problem of extending $g'$ over $D_i$, such that the consequence of (\ref{g'comz}) holds over $D_i$. We can find a complex coordinate $\xi_i$ on $D_i$ such that
$$f^{m_i}(\xi_i) = \exp(2\pi i \frac{n_i}{l_i}) \xi_i.$$
Let $0<k_i<l_i$ be given by $k_in_i =1 \text{ mod } l_i$.
From the assumption that $[\phi_1^{m_i}(p_i)]= (z')^{m_i}[\phi_2^{m_i}(p_i)],$ we conclude that 
$$\left[ \prod_{j=0}^{m_i-1} \varphi_1(f^j(p_i))\right] = \left[ (z')^{m_i} \prod_{j=0}^{m_i-1} \varphi_2(f^j(p_i))\right] .$$
Choose $g_i\in G$ such that 
$$ \prod_{j=0}^{m_i-1} \varphi_1(f^j(p_i)) = (z')^{m_i} g_i^{-1}\prod_{j=0}^{m_i-1} \varphi_2(f^j(p_i)) g_i.$$
Now choose a curve $\gamma_i(t) \in G$, $t\in [0,1]$ connecting $g_i$ and $g'(\xi_i^{-1}(1))$.
Let $I_i = \xi_i^{-1}([0,1])$ and define $g'$ over $I_i$ by the assignment $g'(\xi_i^{-1}(t)) = \gamma_i(t)$. We use (\ref{g'comz}) to extend $g'$ to $f^{jm_i}(I_i)$, $j=0, \ldots, l_i-1$. This makes the extension consistent by the choice of $\gamma_i(0)$. Since $G$ is simply connected, it is possible to
extend $g'$ over $\xi_i^{-1}(\{re^{i\theta}| 0 \leq r \leq1, \mbox{  } 0\leq \theta\leq \frac{2\pi}{l_i}\})$. Now use the equation (\ref{g'comz}) to extend $g'$ to the rest of $D_i$ and to $f^j(D_i)$, $j=1, \ldots m_i$. We have now extended $g'$ to become an element $g\in {\mathcal G}$ which solves the equation $g\phi_1 g^{-1} = z'_0 \phi_2$ globally on $\Sigma$. Thus we have proved that $H$ is injective.

For the surjectivity of the map $H$, suppose we are given $(z, c_1,\ldots c_n)\in Z\times C^n$. Choose $\varphi': \Sigma' \ra G$ to be a constant map, such that $\varphi'^m = z$. We will now extend $\varphi'$ to a map $\varphi : \Sigma \ra G$, such that $\phi(x,g) = (f(x), \varphi(x)g)$ is a lift of $f$ and such that $H(\phi) = (z, c_1,\ldots c_n)\in Z\times C^n$. Choose $\varphi(p_i)\in c_i$, for $i=1, \ldots n$. Further for each $i=1, \ldots n$, choose a curve $g(t)$, $t\in [0,1]$ in $G$ from $\varphi(p_i)$ to $\varphi'$, which is constant near $0$ and $1$. Now define $\varphi(f^j(\xi_i^{-1}(t))) = g(t)$ for $i=1, \ldots n$ and $j=0, \ldots m-2$ and
$\varphi(f^{m-1}(\xi_i^{-1}(t))) = z g(t)^{(1-m)}$ for $i=1, \ldots n$ and for $t\in [0,1]$. Now extend $\varphi$ to $\xi_i^{-1}(\{re^{i\theta}| 0 \leq r \leq1, \mbox{  } 0\leq \theta\leq \frac{2\pi}{l_i}\})$ for $i=1, \ldots n$ in such a way that $\varphi$ is constant near the boundary. This is possible since $G$ is simply connected. Further extend $\varphi$ by defining $\varphi(f^j(p)) = \varphi(p)$ for all $p\in \xi_i^{-1}(\{re^{i\theta}| 0 \leq r \leq1, \mbox{  } 0\leq \theta\leq \frac{2\pi}{l_i}\})$ for $i=1, \ldots n$ and $j=1,\ldots m-2$. Finally, we
set 
$$\varphi(f^{m-1}(p)) = z \varphi(p)^{m-1}$$
for all $p\in \xi_i^{-1}(\{re^{i\theta}| 0 \leq r \leq1, \mbox{  } 0\leq \theta\leq \frac{2\pi}{l_i}\})$ for $i=1, \ldots n$. We have now extended $\varphi$ to all of $\Sigma$, such that it has the wanted properties.

\endproof

From now on we will use $H$ to identify $\Delta$ with $\{(z,c_1, \ldots, c_n) \in Z \times C^n \mid c_i^{l_i} = z\}/Z$.

Let us turn to the description of the subsets $|{\mathcal M}'|_\delta$. Recall that $\epsilon$ denotes a curve in $\Sigma'$ which connects two points related by $f$. It projects to a closed curve in $\tilde \Sigma'$, which represents a generator, say $[\epsilon]$, of $\pi_1(\tilde \Sigma')/\pi_*(\pi_1(\Sigma')) \cong {\mathbb Z}_m$.

Let 
$$Z_{(c_1,\ldots, c_n)} =\{ \rho\in \Hom( \pi_1(\tilde \Sigma')/\pi_*(\pi_1(\Sigma')) , Z) \mid \rho([\epsilon])^{-m_i}c_i = c_i\mbox{, } i = 1, \ldots n\}.$$
Evaluation on $\epsilon$ renders $Z_{(c_1,\ldots, c_n)}$ a subgroup of $Z$. Let further ${\mathcal M}(\tilde{\Sigma}', c_1^{-k_1},\ldots,c_n^{-k_n})$ be the moduli space of flat $G$-connections on $\tilde\Sigma'$ with holonomy around $\partial_i(\tilde{\Sigma}')$ in $c_i^{k_i}$, $i=1,\ldots, n$ and trivial holonomy around $\partial_0(\tilde{\Sigma}')$.  For each $\delta\in \Delta$ we choose an element $(z,c_1, \ldots, c_n) \in \{(z,c_1, \ldots, c_n) \in Z \times C^n \mid c_i^{l_i} = z\}$ representing $H(\delta)$ and we then define
$$Z_\delta = Z_{(c_1,\ldots c_n)}$$
and
$$c(\delta) := (c_1^{-k_1}, \ldots c_n^{-k_n}).$$
Denote by ${\mathcal M}''(\tilde{\Sigma}', c(\delta))$ the space of flat $G$-connections in ${\mathcal M}'(\tilde{\Sigma}', c(\delta))$ which remain irreducible when pulled back via $\pi$.

\begin{theorem}\label{fpi}
Pull back with respect to $\pi$ followed by pull back with respect to a certain broken gauge transformation $g$ constructed below, followed by extension to all of $\Sigma$ induces an isomorphism
$$|{\mathcal M}'|_\delta \cong {\mathcal M}''(\tilde{\Sigma}', c(\delta))/Z_{\delta}.$$
\end{theorem}

\proof

First pick $(z,c_1,\ldots,c_n)$ representing a point in $\Delta$ and further consider a lift $\phi$ such that $H(\phi) = (z,c_1,\ldots,c_n)$.

We consider the restriction map $r : {\mathcal A}'_{\phi} \rightarrow {\mathcal A}'_{\phi}(P')/{\mathcal G}_\phi(P')$ induced by the restriction to $\Sigma'\subset \Sigma$. Suppose we have $A_1,A_2 \in {\mathcal A}'_{\phi}$ such that $r(A_1) = r(A_2)$. That means that $r(A_1)$ and $r(A_2)$ induces the same representation of $\pi_1(\Sigma')$ into $G$ modulo conjugation. But then they induce the same representation of $\pi_1(\Sigma)$ into $G$ modulo conjugation. Hence there exists $g\in {\mathcal G}(P)$ such that $g^*A_1 = A_2.$ From this we see immediately that $g\in {\mathcal G}_\phi.$ Furthermore the restriction map is of course  ${\mathcal G}_\phi$-invariant and we conclude that 
$$r : |{\mathcal M}'|_\delta \ra {\mathcal A}'_{\phi}(P')/{\mathcal G}_\phi(P')$$
is injective.

By the above proof of Proposition \ref{H}, we observe that we can assume that $\phi'= f\times \varphi$ over $\Sigma'$, where $\varphi\in G$ solves the equation $\varphi^m = z$. 
Consider now a map $u: \Sigma' \ra U(1)$ with the property that $u(f(x)) = e^{\frac{2\pi i}{m}} u(x)$ for all $x\in \Sigma'$. Since $\Sigma' \ra \tilde\Sigma'$ is a cyclic cover, with $\langle f \rangle$ as its Galois group, it is easy to construct such a map. We use $u$ to define a broken gauge transformation $g$ of $P'$, which has a discontinuity along $u^{-1}(1)$. First choose a curve $\varphi_t \in G$, $t\in [0,1]$ such that $\varphi_t = 1$ near $t=0$ and $\phi_t = \varphi$ near $t=1$. Now define $h : U(1) \ra G$ by 
$$h(e^{\frac{2\pi i}{m} (j + t)}) = \varphi^{-j}\varphi(t)^{-1}\mbox{, } t\in [0,1].$$
We observe that $h$ has a discontinuity at $1\in U(1),$ where it jumps multiplicatively by $z$. Let $g= h\circ u$. We observe that $g$ acts on ${\mathcal A}(P')$ and it induces a bijection
$$g^* : {\mathcal A}_\phi(P') \ra {\mathcal A}_{\phi'_{\rm  tr}}(P')$$
where $\phi'_{\rm tr} = f|_{\Sigma'}\times \Id$. Also conjugation by $g$ induces a bijection between ${\mathcal G}_{\phi'}(P')$ and ${\mathcal G}_{\phi'_{\rm tr}}(P')$, which is compatible with the identification of ${\mathcal A}_{\phi'}(P')$ and ${\mathcal A}_{\phi'_{\rm tr}}(P')$ via $g^*$. Let $P'_{\phi'_{\rm tr}} = P'/\langle \phi'_{\rm tr} \rangle$, which is clearly a principal $G$-bundle over $\tilde \Sigma'$. Let us moreover denote by ${\mathcal A}''(P'_{\phi'_{\rm tr}})$ the space of flat $G$-connections in $P'_{\phi'_{\rm tr}}$ which remain irreducible when pull back to $\Sigma'$ via $\pi$.
We then get the identification
$$\pi^*: {\mathcal A}''(P'_{\phi'_{\rm tr}}) \stackrel{\cong}{\ra} {\mathcal A}'_{\phi'_{\rm tr}}(P') $$
via pull back to $\Sigma'$. Moreover, we get a short exact sequence of groups
$$ 1 \ra {\mathcal G}(P'_{\phi'_{\rm tr}}) \ra {\mathcal G}_{\phi'_{\rm tr}}(P') \ra \Hom( \pi_1(\bar \Sigma')/\pi_*(\pi_1(\Sigma')) , Z) \ra 1,$$
where $g'\in  {\mathcal G}_{\phi'_{\rm tr}}(P')$ gets maps to the homomorphism $\rho$ given by $\rho([\epsilon]) = g'(\epsilon_1)/g'(\epsilon_0).$

Recalling the notation from the proof of Proposition \ref{H}, we let $q^0_i = \xi_i^{-1}(1)$ and $q_i = f^{m_ik_i}(q_i^0).$

Suppose $A\in {\mathcal A}'_\phi(P')$.
Let $PT_A(-,+) : P_- \rightarrow P_+$ denote parallel transport with respect to
$A$ in $P$ along the curves specified in the following figure. 

\bigskip
\bigskip

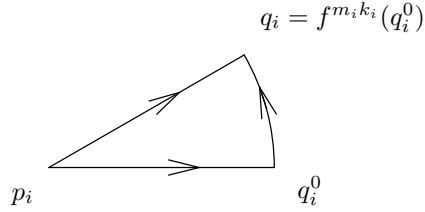
\begin{center}
\begin{texdraw}
\drawdim cm \linewd 0.02
\arrowheadtype t:V 
\avec (2 0) \lvec (3 0)
\move (0 0) \avec (1.732 1) \lvec (2.598 1.5)
\move (0 0) \larc r:3 sd:0 ed:30
\move (2.897 0.776) \avec (2.803 1.0657)
\htext (-.5 -.5) {$p_i$}
\htext (3.3 -.5) {$q_i^0$}
\htext (2.8 1.8) {$q_i= f^{m_ik_i}(q_i^0)$}
\end{texdraw}
\end{center}

\bigskip
\bigskip

By the invariance of $A$ under ${\phi}$, we get the following
commutative diagram \[\begin{array} {ccccccc} & P_{p_i}
&\stackrel{{\phi}^{-m_ik_i}(p_i)}{\leftarrow} & P_{p_i}&
\stackrel{\mbox{Id}}{\leftarrow} & P_{p_i}&\\ PT_A(p_i,q_i^0) & \downarrow
&\mbox{ } \mbox{ } PT_A(p_i,q_i)&\downarrow & &\downarrow&PT_A(p_i,q_i^0)\\ &
P_{q^0_i} &\stackrel{{\phi}^{-m_ik_i}(q_i)}{\leftarrow} & P_{q_i} &
\stackrel{PT_A(q^0_i,q_i)}{\leftarrow} & P_{q^0_i}& \\ \end{array}\] 
hence
\[{\phi}^{-m_ik_i}(q_i) PT_A(q_i^0,q_i)=
PT_A(p_i,q_i^0){\phi}^{-m_ik_i}(p_i)PT_A(p_i,q_i^0)^{-1}.\]

Now let ${\tilde A}$ be the unique connection in $P'_{\phi'_{\rm tr}}$ such that $g^*\pi^*({\tilde A}) = r(A)$. Then by the above formula, we see that 
the holonomy of $\tilde A$ around the $i$'th boundary component of $\tilde \Sigma'$ is in the conjugacy class $c_i^{-k_i}$.

Suppose conversely, we have a connection ${\tilde A}$ in $P'_{\phi'_{\rm tr}}$  which has holonomy around $\partial_i({\tilde \Sigma}')$ in $c_i^{-k_i}$, $i=1,\ldots, n$ and trivial holonomy around $\partial_0({\tilde \Sigma}')$. Then we will construct an $A\in {\mathcal A}_\phi$ such that $g^*\pi^*({\tilde A}) = r(A)$. We consider  $A'=g^*\pi^*({\tilde A})$ and let $h^{(i)}$ be the parallel transport of $A'$ along $\xi_i^{-1}(e^{2\pi i t})$, for $t\in [0,1/l_i]$. Then $\phi^{-m_ik_i}(q_i)h^{(i)} \in c^{-k_i}_i$ by definition. Pick $p^{(i)}\in G$ such that $p^{(i)} \phi^{-m_ik_i}(p_i) (p^{(i)})^{-1} = \phi^{-m_ik_i}(q_i)h^{(i)} $. Choose a curve $p_t^{(i)}$ from $1\in G$ to $p^{(i)}$, which is constant near $t=0$ and near $t=1$. Pick $A$ along $\xi_i^{-1}([0,1])$ such that $p_t^{(i)}$ is covariant constant. Then the parallel transport of $A$ along $\xi_i^{-1}([0,1])$ is $p^{(i)}$. We define $A$ along $\xi_i^{-1}( e^{2\pi i/k_i} [0,1])$ to be $A = (\phi^{-m_ik_i})^*A.$ Then we have the commutative diagram

\[\begin{array} {ccccccc} & P_{p_i}
&\stackrel{{\phi}^{-m_ik_i}(p_i)}{\leftarrow} & P_{p_i}&\\ PT_A(p_i,q_i^0) & \downarrow
&\mbox{ } \mbox{ } PT_A(p_i,q_i)&\downarrow & \\ &
P_{q^0_i} &\stackrel{{\phi}^{-m_ik_i}(q_i)}{\leftarrow} & P_{q_i} & \\ \end{array}\] 

Now we compute
\begin{eqnarray*}
PT_A(p_i,q_i) PT_A(p_i,q^0_i) ^{-1} & = & \phi^{m_ik_i}(q_i) PT_A(p_i,q^0_i)  \phi^{-m_ik_i}(p_i)PT_A(p_i,q^0_i) ^{-1}\\ & = & PT_A(q^0_i,q_i).
\end{eqnarray*}
But then we can simply extend $A$ over the sector $(p_i,q_i^0,q_i)$ such that $A$ is trivial near $p_i$ and such that near the boundary, it is a pull back of $A$, from the boundary of the sector with respect to some smooth map onto the boundary. Now we extend $A$ to all the remaining sectors near $p_i,f(p_i), \ldots f^{m_i-1}(p_i)$ by requiring  invariance under $\phi$. Doing this for all $i=1, \ldots, n$ and on $\pi^{-1}(D_0)$ in a $\phi$-invariant way, gives $A\in {\mathcal A}_\phi$ such that $g^*\pi^*({\tilde A}) = r(A)$ as required. The Theorem now follows by observing that the subgroup of ${\mathcal G}_{\phi_{\rm tr}}(P')$ which maps $r({\mathcal A}_\phi)$ to itself consists exactly of the gauge transformations $g'\in {\mathcal G}_{\phi_{\rm tr}}(P')$ such that
$(g'(\epsilon_1)/g'(\epsilon_0))^{m_ik_i}c_i^{k_i} = c_i^{k_i}.$
\endproof

Let us use the notation $|{\mathcal M}|_{\delta}$ for the subset of  
$|{\mathcal M}|$ which can be represented by a lift which represents 
$\delta\in \Delta$ and $\check{\pi} : {\mathcal M}''(\tilde{\Sigma}', c(\delta))/{Z_\delta} \rightarrow |{\mathcal M}|_\delta$ for the map constructed in the proof of Theorem \ref{fpi}.

\begin{prop}  \label{ici}
We have that the map
$$ \check{\pi} :  {\mathcal M}(\tilde{\Sigma}', c(\delta))/Z_{\delta}\rightarrow |{\mathcal M}|_\delta.$$
is surjective and injective on the subset ${\mathcal M}''(\tilde{\Sigma}', c(\delta))/Z_{\delta}$. In particular if $|{\mathcal M}|_\delta \subset {\mathcal M}'$ then $\check{\pi}$ is an isomorphism onto the component $|{\mathcal M}|_\delta$.
\end{prop}

\begin{remark}  \label{Recomplex} We notice that the map $\check{\pi}$ is  holomorphic with respect to the complex structure induced on $|{\mathcal M}|_\delta$ from $\sigma$ and the complex structure induced on ${\mathcal M}(\tilde{\Sigma}', c(\delta))$ by $\tilde{\Sigma}_\sigma \cong \Sigma_\sigma/\langle f \rangle$. The moduli space ${\mathcal M}(\tilde{\Sigma}', c(\delta))$ with this complex structure is denote ${\mathcal M}_\sigma(\tilde{\Sigma}', c(\delta))$.
\end{remark}

\begin{remark}  \label{Resymp} Let
$\omega(c(\delta))$ be the natural symplectic structure on 
${\mathcal M}'({\tilde \Sigma}',c(\delta))$ and $\omega_{\delta} = \omega|_
{|{\mathcal
M}'|_{\delta}}$. It is easy to verify that 
\[\omega_{\delta} = m \check{\pi}^*(\omega(c(\delta))).\]
\end{remark}

\begin{remark}  In the special case where $f$ has no special orbits, i.e. 
$\Sigma\rightarrow {\tilde \Sigma}$ is a principal ${\Bbb Z}_m$ bundle, 
we get that
\[\Delta = Z/Z,\]
where $z'\in Z$ acts on $Z$ by $z\mapsto (z')^m z$ and for each $\delta \in \Delta$, we have that $Z_\delta = Z_m$, where
$$ Z_m = \{z\in Z \mid z^m =1\}.$$
Further $$ \check{\pi} : {\mathcal M}(\tilde{\Sigma})/Z_m \rightarrow |{\mathcal M}|_\delta$$
is a surjection.

\end{remark}

Let us now concentrate on the part of the fixed point set of $f$ which is 
contained in the reducible locus of ${\mathcal M}$.
Let $A$ be a connection representing a point in $|{\mathcal M}|$. From our earlier 
discussion, we know there is a lift ${\phi}$ of $f$ such that ${\phi}A=A$ 
and ${\phi}^m\in Z_A.$ Note that any other such lift is of the form 
${\phi}g$ where $g\in Z_A$.
Let $C$ be the set of connected components of $|{\mathcal M}|$ and for each $c\in C$ let $|{\mathcal M}|^c$ be the corresponding component subset of $|{\mathcal M}|$. We note of course 
that there is a map from $\Delta$ to $C$, by mapping $\delta\in \Delta$ to 
the 
component of $|{\mathcal M}|$ which contains $|{\mathcal M}'|_\delta.$ It induces an equivalence relation $\sim$ defined by $\delta_1 \sim \delta_2$ if and only if $|{\mathcal M}'|_{\delta_1}$ and  $|{\mathcal M}'|_{\delta_2}$ are contained in the same component of $|{\mathcal M}|$.
We can in general not expect that $\sim$ is trivial, since an element $g\in Z_A$ 
may map a lift representing say $\delta_1\in \Delta$ to another lift which 
represents 
a different element say $\delta_2\in \Delta$, which would mean that 
$|{\mathcal M}|_{\delta_1}$ and $ |{\mathcal M}|_{\delta_2}$ intersect in the 
reducible 
locus of ${\mathcal M}$.
This reflects the fact that $\Delta$ indexes some of the irreducible components of  
$|{\mathcal M}|$, whereas $C$ indexes the connected components.

Observe that the map from $\Delta$ to $C$ is surjective iff $|{\mathcal M}|=\overline{|{\mathcal M}'|}$.
 In order to establish that $\bigcup_{\delta\in\Delta} |{\mathcal M}|_\delta$  equals $|{\mathcal M}|$, we just 
need to show the following:
For any $A\in {\mathcal A}_P$ and for any lift $\phi$ of $f$, such that $\phi^m \in Z_A$, we can choose $g\in Z_A$ such that $({\phi}g)^m\in Z.$ -- The map 
$Z_A \rightarrow Z_A$ given by $g \rightarrow ({\phi}g)^m$ is a covering map, 
since its derivative is a composition of isomorphisms. To ensure that we can find 
a $g\in Z_A$ such that  $({\phi}g)^m\in Z$, it suffices to prove that $Z_A/Z$ 
is connected. Since the stabilizer of $A$ is the centralizer of the holonomy 
group of $A$, we obtain the desired description of $|{\mathcal M}|$ for $G=SU(N)$
by proving the following Proposition.

\begin{prop} \label{conn}
The group $Z(H)/Z$ is connected for all Lie-subgroups $H\subset SU(N).$
\end{prop}

\begin{remark} If we replace $SU(N)$ by say $SO(N)$ this Proposition is 
false.
\end{remark}

\proof

We think of the inclusion $\rho : H \rightarrow SU(N) \subset U(N)$ as a 
representation of $H$. We know that we can write $\rho$ as a direct sum of 
irreducible representations
\[\rho = \bigoplus_{i=1}^r \rho_i.\]
The automorphism group of the representation $\rho$, as a unitary 
representation is now very easy to describe. The automorphism group of $\rho_i$ is by Schur's Lemma ${\Bbb C}^*.$ Assume 
that we have ordered the representations such that there exists $i_1,\ldots,i_u 
\in \{1,\ldots,r\}:$
\[\rho_i = \rho_j \mbox{ iff } \exists l :  i_l \leq i,j <i_{l+1}.\]
Let $n_j = i_{j+1} - i_j,$ and $m_j = \dim \rho_j.$ Note that $N=\sum_
{j=1}^u n_jm_j.$

An automorphism of $\rho$ has to respect the decomposition of $\rho$.
From this we see that the centralizer of $\rho(H)$ in $U(N)$ is isomorphic 
to
\[\prod_{j=1}^u U(n_j).\]
Hence we get the following expression
\[Z(H) \cong d^{-1}(1),\]
where
\[d : \prod_{j=1}^u U(n_j) \rightarrow U(1),\]
is given by
\[d(A_1,\ldots,A_u) = \prod_{i=1}^u \det(A_i)^{m_j}.\]

Notice that d is a principal $Z(H)$ bundle. Therefore we also have that 
$\prod_{j=1}^u U(n_j)/Z \rightarrow U(1)$ is a principal $Z(H)/Z$-bundle. 
Clearly $\prod_{j=1}^u U(n_j)/Z$ is connected, so by the long exact sequence 
of homotopy groups for the principal bundle $\prod_{j=1}^u U(n_j)/Z 
\rightarrow U(1)$ we just need to lift a representative of a generator of
$\pi_1(U(1))$ to a curve in $\prod_{j=1}^u U(n_j)/Z$ such that it starts in 
$1\in\prod_{j=1}^u U(n_j)/Z$ and ends in $Z\subset \prod_{j=1}^u U(n_j)/Z$. 
If $A_i(t) = \expo(2\pi i \frac{t}{n})I$ then
\[d(A_1(t),\ldots,A_u(t)) = \prod_{i=1}^u \expo(2\pi i \frac{n_i m_j}{n}t) 
=\expo(2\pi i t),\]
so $(A_1(t),\ldots,A_u(t))$ is exactly such a curve, hence the Proposition is 
proved.

\begin{flushright} $\Box$ \end{flushright}

We collect the result for $G=SU(N)$ in the following
\begin{theorem}\label{TSUN}

For $G=SU(N)$ we have the following decomposition of $|{\mathcal M}|$ into
irreducible components:
\[|{\mathcal M}| = \bigcup_{\delta\in\Delta}|{\mathcal M}|_\delta,\]
where $\Delta= \{(z,c_1, \ldots c_n) \in Z \times C^n\mid c_i^{l_i} = z\}/Z,$  
\[\check{\pi} : {\mathcal M}(\tilde{\Sigma}', c(\delta))/Z_{\delta}\rightarrow |{\mathcal M}|_\delta,\]
is the surjection defined just above Proposition \ref{ici} and $c(\delta)$ and $Z_\delta$ are defined just above Theorem \ref{fpi}. This map is holomorphic with respect to the natural
complex structures induced by $\sigma$. The map $\check{\pi}$ is injective on the subset ${\mathcal M}''(\tilde{\Sigma}', c(\delta))/Z_{\delta}$. In particular if $|{\mathcal M}|_\delta \subset {\mathcal M}'$ then $\check{\pi}$ is an isomorphism onto the component $|{\mathcal M}|_\delta$.
\end{theorem}

In \cite{AGr} J. Grove, joint with this author, investigated the fixed point set for an automorphism of a Riemann surface acting on the moduli space 
of rank two semi-stable holomorphic bundles with fixed determinant. If we choose $G=SU(2)$, we see that the above description of the fixed point set  is equivalent to the description of the fixed points given in Theorem 3.4 in \cite{AGr}. In that paper we further describe the fibers of $\check{\pi}$ over $|{\mathcal M}|_\delta - |{\mathcal M}'|_\delta$ in the rank two case using algebraic geometric techniques.

\section{The moduli space of flat connections on a mapping torus.}

\label{3MF} In this Section we shall examine the relation between the moduli
space of flat $G$-connections on $\Sigma_f$, denoted ${\mathcal M}(\Sigma_f)$ and
the fixed point set $|{\mathcal M}|$. In the first part of this Section we have
no assumptions on the diffeomorphism $f$, except it has to be orientation
preserving. It is obvious from the construction of $\Sigma_f$ that
\[\pi_1(\Sigma_f) = \langle \pi_1(\Sigma), \gamma \mid \forall \alpha\in \pi_1(\Sigma) :
\gamma^{-1} \alpha \gamma= f\alpha\rangle.\]

There is a natural map \[r : {\mathcal M}(\Sigma_f) \rightarrow |{\mathcal
M}|,\] given by restricting a flat connection on $\Sigma_f$ to $\Sigma \times
\{0\}$. Recall that ${\mathcal M}(\Sigma_f)^c =
r^{-1}(|{\mathcal M}|^c),$ hence ${\mathcal M}(\Sigma_f)^c$ is a union of
components which map to the component $|{\mathcal M}|^c$.

Suppose $[\rho]\in |{\mathcal M}|$. We then have that there exists $g\in G$
such that $g^{-1}\rho g = f(\rho)$. Note that $g\in N(\rho(\pi_1(\Sigma)))$. We want to compute
$r^{-1}([\rho])$. Suppose $[\tilde{\rho}]\in r^{-1}([\rho])$. We may assume that
$\tilde{\rho}$ is such that $\tilde{\rho}|_{\pi_1} =\rho$. We then have that
\[\tilde{\rho}(\gamma) \rho \tilde{\rho}(\gamma)^{-1} = f(\rho),\] so
\[\tilde{\rho}(\gamma) g^{-1} \in Z_\rho.\] Hence we see that
$\tilde{\rho}(\gamma)\in Z_\rho g$ parametrizes the different homomorphisms
$\tilde{\rho} : \pi_1(\Sigma_f) \rightarrow G$ such that $\tilde{\rho}|_{\pi_1}
=\rho$.

Suppose that $\tilde{\rho}_i$ are two such representations of $\pi_1(\Sigma_f)$
which are conjugate, say $\tilde{g}\tilde{\rho}_1\tilde{g}^{-1} =
\tilde{\rho}_2.$ Then $\tilde{g}\in Z_{\rho}$. If $\tilde{\rho}_1(\gamma) = z
g,$ then we have that $\tilde{\rho}_2(\gamma)= \tilde{g} z g \tilde{g}^{-1}$.
Since $g\tilde{g}^{-1} g^{-1}\in Z_\rho$ we see that $\tilde{\rho}_2(\gamma)\in
Z_\rho g$.

We conclude that \[r^{-1}([\rho])\cong Z_\rho g / Z_\rho,\] where $Z_\rho$ acts
on $Z_\rho g$ by conjugation. In the generic case, when $Z_\rho= Z(G)$, we thus
see that $r^{-1}([\rho]) \cong Z(G).$ So from this we see that $r$ is
a $|Z(G)|$-sheeted cover over $|{\mathcal M}'|$, but it has singularities over the reducible part of
$|{\mathcal M}|$.

From this we see that if $|{\mathcal M}|^c$ is a component of $|{\mathcal
M}|$ contained in ${\mathcal M}'$, then $r^{-1}(|{\mathcal M}|^c)$ is a union of
smooth components of ${\mathcal M}(\Sigma_f)$ which together form a
$|Z(G)|$-fold cover of $|{\mathcal M}|^c$.

Let $C_f$ be the set of connected components of ${\mathcal M}(\Sigma_f)$. 
The map $r$ induces a map $r : C_f \rightarrow C$. For $SU(N)$, we have seen 
that $Z_\rho/Z$ is connected (see Prop. \ref{conn}), so $|r^{-1}(c)|\leq |Z|$ for all $c\in C$.

Let $[\rho]\in|{\mathcal M}|.$ From the discussion in the beginning of the
previous Section, we know that there is a lift ${\phi}$ of $f$ to $P$ and a
connection $A$ in $P$ invariant under ${\phi}$ such that $[A]=[\rho]$. Let
$P_{{\phi}}$ be the mapping torus of ${\phi}$. Then $A$ induces naturally
a flat connection in $P_{{\phi}}$, which we denote $A_{{\phi}}$. By
standard theory any flat connection in a principal G-bundle over $\Sigma_f$ is
equivalent to some $A_{{\phi}}$ for some A and some ${\phi}$. Hence we can
represent all $[\tilde{\rho}]\in r^{-1}([\rho])$ this way.

Recall that we have a natural action of $f$ on ${\mathcal L}$ covering $f$'s 
action on ${\mathcal M}$. Since $[A]\in|{\mathcal M}|$, we see that 
$f:{\mathcal L}^k_{[A]} \rightarrow {\mathcal L}^k_{[A]}$. The following 
Lemma gives a formula for the trace of $f$'s action on ${\mathcal L}^k_{[A]}$.

\begin{lemma}\label{TCS}
 We have that
 \[\tr(f : {\mathcal L}^k([A])\rightarrow {\mathcal
L}^k([A])) = \expo(2\pi i k \CS(P_{{\phi}}, A_{{\phi}})).\] 
\end{lemma}

\begin{remark}
From the formula, we see that $\CS(P_{{\phi}}, A_{{\phi}})$ mod 
${\Bbb Z}$ only depends on $r(A_{{\phi}})$ and therefore is constant on 
$r^{-1}(|{\mathcal M}|^c)$. We therefore use the notation
\[\CS(\Sigma_f,c) = CS(P_{{\phi}},A_{{\phi}})\]
whenever $r(A_{{\phi}})\in |{\mathcal M}|^c$.
\end{remark}

\proof
This formula follows directly from Theorem 2.19 in \cite{F}.

\begin{flushright} $\Box$ \end{flushright}
Using a standard Mayer-Vietoris argument, we obtain the formula
\begin{equation} \dim_{\Bbb R} E_1({\phi},1) =
\dim_{\Bbb R} H^1(\Sigma_f,d_{A_{{\phi}}}) - \dim_{\Bbb R}
H^0(\Sigma_f,d_{A_{{\phi}}}), \label{dFC} \end{equation}
where $E_1({\phi},1)$ is the 1-eigenspace of ${\phi}$ on
$H^1(\Sigma,d_A).$ 
Assume again that $f$ is of
finite order. We get this way a formula for the local dimension of $|{\mathcal
M}|$ around $[\rho]$ which of course is equal to $\dim_{\Bbb R}
E_1({\phi},1)$. Let us record this observation in the following Lemma.
\begin{lemma} \label{dc} Let 
\begin{equation}
d_c = \max_{[A]\in {\mathcal M}(\Sigma_f)^c} \frac{1}{2}
(\dim_{\Bbb R} H^1(\Sigma_f,d_{A}) - \dim_{\Bbb R} H^0(\Sigma_f,d_{A})),
\label{dF}
\end{equation}
where $\max$ means the generic maximum, i.e. the maximum value taken on a
Zariski open subset of ${\mathcal M}(\Sigma_f)^c$. If $f$ is of finite order, we
then have that 
\[ d_c = \dim_{\Bbb C} |{\mathcal M}(\Sigma)|^c.\]

\end{lemma}

\section{Localization of the Witten invariant.}

\label{L} 
In this Section we are going to apply the Lefschetz-Riemann-Roch
formula for singular varieties due to Baum, Fulton, McPherson and Quart to 
give an expression for the Witten invariant of
any finite order mapping torus. See \cite{BFM} and \cite{BFQ} for the general
statement of their Theorem.

In order to apply the Lefschetz-Riemann-Roch formula, we have to determine the 
trace of $f$'s action on the higher
cohomology groups of ${\mathcal L}^k_\sigma$ over ${\mathcal M}_\sigma$. The following 
well-known 
Theorem takes care of that problem.

\begin{theorem} 
 \[H^i({\mathcal M}_\sigma,{\mathcal L}^k_\sigma) = 0\mbox{ for } i>0.\] 
\end{theorem}

\proof
This fact follows from Kodaira vanishing which is applicable since ${\mathcal M}_\sigma$
 has rational singularities and the fact that the canonical bundle is negative
(Theorem F. 
in \cite{DN}).
\begin{flushright} $\Box$ \end{flushright}

\begin{theorem}\label{TWF}
 We have the following expression for the Witten invariant of the
mapping torus $\Sigma_f$ calculated with respect to the Atiyah 2-framing of
$\Sigma_f$. 
\begin{eqnarray} Z^{(k)}_G(\Sigma_f) & = &\Det(f)^{-\frac{1}{2}\zeta}
\sum_{c\in C} \expo(2 \pi i k\CS(\Sigma_f,c)) \expo(k\omega_c) \cap
\tau_\cdot(L_\cdot^c({\mathcal O}_{{\mathcal M}_\sigma}))\nonumber\\ & = &
\Det(f)^{-\frac{1}{2}\zeta} \sum_{c\in C} \expo(2 \pi i k\CS(\Sigma_f,c))
\sum_{i=0}^{d_c} \frac{1}{i!}((\omega_c)^i \cap \tau_i(L_\cdot^c({\mathcal
O}_{{\mathcal M}_\sigma}))) k^i \label{LF} 
\end{eqnarray} 
Here $\omega_c$ is the
restriction of $c_1({\mathcal L})$ to $|{\mathcal M}_\sigma|^c$, 
$\tau_\cdot( L_\cdot^c({\mathcal O}_{{\mathcal M}_\sigma}))\in H_\cdot(|{\mathcal M}_\sigma|^c,{\Bbb C})$ 
is
 the homology class defined in the Theorem on page 180 in \cite{BFM} and 
in \S 2 in \cite{BFQ}, $d_c$ is given by (\ref{dF}), $C$ is the set of 
connected 
components of $|{\mathcal M}_\sigma|$, which was described explicitly in chapter 
\ref{FPS}.
 The framing 
correction term $\Det(f)^{-\frac{1}{2}\zeta}$ is given by (\ref{FCT}).

\end{theorem}

\proof
Recall that we are using the following notation for the component decomposition
of the fixed point set of $f$ on ${\mathcal M}_\sigma$
\[|{\mathcal M}_\sigma| = \bigcup_{c\in C} |{\mathcal M}_\sigma|^c.\]
Following the notation of \cite{BFM} and \cite{BFQ} we denote the Grothendieck group of all equivariant coherent sheaves on ${\mathcal M}_\sigma$ by
$K^{\eq}_0({\mathcal M}_\sigma)$ and Grothendieck group of all coherent sheaves on $|{\mathcal M}_\sigma|^c$ by $K^{\alg}_0(|{\mathcal M}_\sigma|^c)$. Let 
\[L_\cdot^c : K^{\eq}_0({\mathcal M}_\sigma) \rightarrow K^{\alg}_0(|{\mathcal
M}_\sigma|^c)\otimes {\Bbb C}\]
be the localizing homomorphism defined in \S 2 in \cite{BFQ}, and 
\[\tau_\cdot : K^{\alg}_0(|{\mathcal M}_\sigma|^c) \rightarrow H_\cdot(|{\mathcal
M}_\sigma|^c)\]
the homomorphism defined in the Theorem on page 180 in \cite{BFM}. The Lefschetz-Riemann-Roch
formula of Baum, Fulton, McPherson and Quart then states that
\[\tr(f : H^0({\mathcal M}_\sigma, {\mathcal L}_\sigma^k) \rightarrow 
H^0({\mathcal M}_\sigma, {\mathcal L}_\sigma^k)) = 
\sum_{c\in C} a_c^k \ch({\mathcal L}^k_\sigma|_{|{\mathcal M}_\sigma|^c}) \cap
\tau_\cdot(L_\cdot^c({\mathcal O}_{{\mathcal M}_\sigma}))\]
where $a_c$ is the
complex number by which $f$ acts on ${\mathcal L}_\sigma|_{|{\mathcal M}_\sigma|^c}.$
As we computed in Lemma \ref{TCS},
\[ a_c = \expo(2\pi i \CS(\Sigma_f,c)).\]
Clearly
\[\ch({\mathcal L}_\sigma^k|_{|{\mathcal M}_\sigma|^c}) 
= \expo(k\omega_c).\]

Hence we arrive at the expression (\ref{LF}).

\begin{flushright} $\Box$ \end{flushright}

Let us now study the perturbative contribution from the smooth components of the
fixed point set. A component $|{\mathcal M}_\sigma|^c $ which is contained in the irreducible part of 
${\mathcal M}_\sigma$ is
smooth and  $|{\mathcal M}_\sigma|^c = |{\mathcal M}_\sigma|_ \delta$ for a unique $\delta\in \Delta$. We get (see Theorem 2.2, 5.  in \cite{BFQ})
\begin{eqnarray} P_c(k) & &\\ \nonumber
& = & \expo(k\omega_c) \cap
\tau_\cdot(L_\cdot^c({\mathcal O}_{{\mathcal M}_\sigma}))k^{-d_c} \\ \nonumber
& = & ( \expo(k
\omega_c) \cup \Ch^\cdot(\lambda_{-1}^c{\mathcal M}_\sigma)^{-1} \cup \Td(T_c)) \cap 
[|{\mathcal M}_\sigma|^c]k^{-d_c}. \label{csc}
\end{eqnarray}
where we have used the notation $T_c = T_{|{\mathcal M}_\sigma|^c}$ for the holomorphic tangent bundle of $|{\mathcal M}_\sigma|^c$ and 
$(\lambda_{-1}^c{\mathcal M}_\sigma)^{-1}\in K_{\alg}^0(|{\mathcal M}_\sigma|^c)\otimes 
{\Bbb C}$ is defined in \S 0.5 in \cite{BFQ}. Let us briefly recall the definition here. Denote by ${\mathcal N}$ the holomorphic conormal bundle to $|{\mathcal M}_\sigma|^c \subset {\mathcal M}_\sigma'$, i.e.
$${\mathcal N} = \left( i_c^*(T{\mathcal M}_\sigma')/T_c\right)^*,$$
where $i_c : |{\mathcal M}_\sigma|^c \ra {\mathcal M}_\sigma'$ is the inclusion map. Let $\lambda_{-1}({\mathcal N}) = \sum (-1)^i \Lambda^i({\mathcal N})$. Then $\lambda_{-1}({\mathcal N})$ determines an element in $K^0_{{\mathbb Z}_m}(|{\mathcal M}_\sigma|^c)$, the Grothendieck group of ${\mathbb Z}_m$-linearized locally free sheaves on $|{\mathcal M}_\sigma|^c$. There is a map from $K^0_{{\mathbb Z}_m}(|{\mathcal M}_\sigma|^c)$ to $K^0_{\alg}(|{\mathcal M}_\sigma|^c)\otimes {\mathbb Z}[{\mathbb C}]$, which is induced by mapping a ${\mathbb Z}_m$-linearized locally free sheaf ${\mathcal F}$ to $\sum {\mathcal F}_a \otimes a$, where ${\mathcal F}_a$ is the eigen-subsheaf of ${\mathcal F}$ corresponding to the eigenvalue $a$. By composing this map with the standard ${\mathbb C}$-value trace on ${\mathbb Z}[{\mathbb C}]$, we get a map to  $K^0_{\alg}(|{\mathcal M}_\sigma|^c)\otimes {\mathbb C}$. The element $(\lambda_{-1}^c{\mathcal M}_\sigma)^{-1}\in K_{\alg}^0(|{\mathcal M}_\sigma|^c)\otimes  {\mathbb C}$ is by definition the inverse of the image of $\lambda_{-1}({\mathcal N}) $ under this composition of maps.

Since $|{\mathcal M}_\sigma|^c = |{\mathcal M}_\sigma|_ \delta\subset {\mathcal M}'$, we get from Proposition \ref{ici} that 
\[|{\mathcal M}_\sigma|^c \cong {\mathcal M}_\sigma(\tilde{\Sigma}', c(\delta))/Z_\delta\] 
and $\omega_c$ gets identified with 
the symplectic structure $m\omega(c(\delta))$ on ${\mathcal M}_\sigma(\tilde{\Sigma}', c(\delta))/Z_\delta$ by Remark \ref{Resymp}. We will in this case use the notation $\omega^c=\omega(c(\delta))$.

We will now express $\Ch^\cdot(\lambda_{-1}^c{\mathcal M}_\sigma)^{-1}$ in 
terms
of known generators for the cohomology of ${\mathcal M}_\sigma(\tilde{\Sigma}', c(\delta))$. Since it is known how to express $\Td(T_c)$ in terms of these 
generators, we get from this an explicit formula for
 the contribution
to the Witten-Reshetikhin-Turaev invariant from ${\mathcal M}_\sigma(\Sigma_f)^c$ expressed in terms of an evaluation of combinations of the known generators for the cohomology ring 
of ${\mathcal M}_\sigma(\tilde{\Sigma}', c(\delta))$.

Let $H_c \subset \Hom(\pi_1(\Sigma), G)$ be the subset of homomorphisms 
 which projects to $|{\mathcal M}_\sigma|^c$. We will now completely explicitly construct a 
universal adjoint bundle $E$ over $|{\mathcal M}_\sigma|^c
\times \Sigma$. Choose a universal cover $\Pi : D \rightarrow \Sigma$ and a lift $f_D : D \ra D$ of $f$.
Now define an action of $\pi_1(\Sigma)\times \overline{G} \times {\mathbb Z}_m$ on $H_c\times D \times \mathfrak{g}_{\mathbb C}$ as follows
$$(\gamma, \overline{g}, i)(\rho,x,X) =(  \overline{g} h_\rho^{-i} \rho h_\rho^{i}\overline{g}^{-1}, f_D^i(x)\gamma, \overline{g}h_\rho^{-i}\rho(\gamma)^{-1}X),$$
for all $(\rho,x,X) \in H_c\times D \times \mathfrak{g}_{\mathbb C}$, where $h_\rho \in \overline{G}$ is uniquely determined by $\rho\circ f_* = h_\rho^{-1} \rho h_\rho$.
Then
\[E = H_c\times D \times \mathfrak{g}/\pi_1(\Sigma)\times \overline{G}\]
is the universal adjoint bundle over $|{\mathcal M}_\sigma|^c
\times \Sigma$ equipped with an action of ${\mathbb Z}_m$ covering the action of $\Id\times \langle f\rangle$ on $|{\mathcal M}_\sigma|^c
\times \Sigma$. Using the complex structure $\sigma$ on $\Sigma$, we naturally get induced the structure of a holomorphic bundle on $E$.

Let $p$ be the projection 
\[p : |{\mathcal M}_\sigma|^c\times \Sigma \rightarrow
|{\mathcal M}_\sigma|^c.\] 
We have that \[p_*(E) = - i_c^*(T_{{\mathcal M}_\sigma}) \in
K^0_{{\Bbb Z}_m}(|{\mathcal M}_\sigma|^c).\] 

According to the Nielsen localization Theorem (see Theorem 4.5 of \cite{Andreas},
which is the analogue of Quart's localization Theorem in \cite{Quart}, but for the
$K$-theory $K^0_{{\Bbb Z}_m}(\cdot)$), the following
diagram is commutative 
\[ \begin{CD} K^0_{{\Bbb Z}_m}(|{\mathcal M}_\sigma|^c\times\Sigma)\otimes_{R({\Bbb Z}_m)}
 \Lambda @>{L_\cdot}>> K_{\alg}^0(|{\mathcal M}_\sigma|^c\times|\Sigma|)\otimes\Lambda\\
 @V{p_*}VV @VV{p_*}V\\ K^0_{{\Bbb Z}_m}(|{\mathcal M}_\sigma|^c)\otimes_{R({\Bbb Z}_m)}
 \Lambda 
@>{L_\cdot}>{\cong}> K^0_{\alg}(|{\mathcal M}_\sigma|^c)\otimes\Lambda. \end{CD}\]
Here $L_\cdot$ is the localization morphism as defined in \cite{Andreas}, $R({\Bbb Z}_m)$ is the representation ring of ${\Bbb Z}_m$, which of course
is $R({\Bbb Z}_m)= {\Bbb Z}[x]/(x^m-1)$ and $\Lambda$ is
the localization of $R({\Bbb Z}_m)$ at the multiplicative set generated by
$(1-x^i), i = 1,\ldots,m-1$ and $|\Sigma|= \{p_1,\ldots,p_l\}$ is the fixed point set of $f$ on $\Sigma$.
Therefore we get the following equality in $ K^{\alg}_0(|{\mathcal
M}_\sigma|^c)\otimes\Lambda$: 
\begin{eqnarray*} L_\cdot(i_c^*(T_{{\mathcal M}_\sigma})) & = & -
\sum_{p\in|\Sigma|} L_\cdot(E_p)\otimes \lambda_{-1}(T_p^*\Sigma)^{-1}\\
& = & - \sum_{s=1}^l L_\cdot(E_{p_s}) (1-x^{n_s})^{-1}. \end{eqnarray*}

From the explicit construction of $E$ as an $f$-linearized bundle we get for 
each fixed point $p_s\in \Sigma$ of $f$ a specific automorphism of $E_{p_s}$ 
whose $m$'th power is the identity. We have that $E_{p_s}^*\in 
K^0_{{\Bbb Z}_m}(|{\mathcal M}_\sigma|^c)$ so
\[L_\cdot(E_{p_s}^*) = \sum_{\nu=0}^{m-1} E_{p_s}^\nu\otimes x^\nu,\]
where $E_{p_s}^\nu$ is the $e^{2\pi i \frac{\nu}{m}}$-eigen subbundle of $E_{p_s}^*$.

Note that $i_c^*(T_{{\mathcal M}_\sigma}) = T_c\oplus {\mathcal N}^*$ in 
$K^0_{{\Bbb Z}_m}(|{\mathcal M}_\sigma|^c)$.

\begin{lemma}\label{calN}
We have the 
following formula for ${\mathcal N}$ in $K^0_{{\Bbb Z}_m}(|{\mathcal M}_\sigma|^c)\cong 
K_{\alg}^0(|{\mathcal M}_\sigma|^c)\otimes \R({\Bbb Z}_m)$
\[
{\mathcal N}  =  \sum_{j=1}^{m-1}  T_c^*\otimes x^j - \sum_{j=0}^{m-1} \left\{ \frac{1}{m}  
\sum_{s=1}^l\sum_{\nu=0}^{m-1} E_{p_s}^\nu\otimes (\mu_m^{-\nu}(n_s) - 
\mu_m^{j-\nu}(n_s)) \right\} \otimes x^j.
\]

\end{lemma}

\proof

Consider the following homomorphism
\[\tr_i : R({\Bbb Z}_m) \rightarrow {\Bbb C}\mbox{, } i=0,\ldots, m-1\]
given by

\[\tr_i(\sum_{j=0}^{m-1}n_j x^j) = \sum_{j=0}^{m-1} n_j e^{2\pi i \frac{ij}{m}}.\]
Let $\varphi : R({\Bbb Z}_m) \rightarrow \Lambda$ be the universal map. Then by universality $\tr_i : \Lambda \rightarrow {\Bbb C}$ is well-defined for $i=1,\ldots,m-1$ such that $\tr_i(p) = \tr_i(\varphi(p))$ for all $p\in R({\Bbb Z}_m)$.
It is easy to verify the following formula for any $p\in R({\Bbb Z}_m)$
\begin{eqnarray*} p & = & \sum_{j=0}^{m-1}\left( \frac{1}{m} \sum_{i=0}^{m-1}
\tr_i(p)e^{-2\pi i \frac{ij}{m}} \right) x^j\\
& = & \sum_{j=0}^{m-1}\frac{1}{m}\left( \tr_0(p) +  \sum_{i=1}^{m-1}\tr_i(p)
e^{-2\pi i \frac{ij}{m}} \right) x^j.
\end{eqnarray*}
Now ${\mathcal N}\in K^0_{{\Bbb Z}_m}(|{\mathcal M}_\sigma|^c)\cong 
K_{\alg}^0(|{\mathcal M}_\sigma|^c)\otimes R({\Bbb Z}_m)$ and say ${\mathcal N} =
\sum_{i=1}^{m-1}{\mathcal N}_i,$ where ${\mathcal N}_i$ is the $e^{2\pi i \frac{i}m}$-eigen-subbundle of ${\mathcal N}$. So it is trivial to verify that
\[{\mathcal N} = \sum_{j=0}^{m-1} \left(
\frac{1}{m}\sum_{i=1}^{m-1}\tr_i(L_\cdot({\mathcal N}))(e^{-2\pi i 
\frac{ij}{m}}-1) \right) \otimes x^j,\]
where $\tr_i : K_{\alg}^0(|{\mathcal M}_\sigma|^c)\otimes R({\Bbb Z}_m) \rightarrow K_{\alg}^0(|{\mathcal M}_\sigma|^c)
\otimes{\Bbb C}$ is just $1\otimes \tr_i$.
Substituting our expression for $L^\cdot({\mathcal N})$ into this formula, we
arrive at the formula stated in the Lemma.

\begin{flushright} $\Box$ \end{flushright}

We have the following general facts.
\[K^0_{{\Bbb Z}_m}(|{\mathcal M}_\sigma|^c) \cong R({\Bbb Z}_m)\oplus
\tilde{K}^0_{{\Bbb Z}_m}(|{\mathcal M}_\sigma|^c) \] 
and the augmentation ideal
$\tilde{K}^0_{{\Bbb Z}_m}(|{\mathcal M}_\sigma|^c)$ is nilpotent. Hence it is clear, that
the units in $K^0_{{\Bbb Z}_m}(|{\mathcal M}_\sigma|^c)$ are precisely, the elements $a$ 
whose
augmentation $\rank(a)$ is invertible in $R({\Bbb Z}_m)$. Let us use the
following notation $\tilde{a} = a -\rank(a)$. For all $a$, such that $\rank(a)$ is invertible, we have that
\[\lambda_{t}(a)^{-1}  =  \lambda_{t}(\rank(a))^{-1}s_{-t}(\tilde{a}),\] 
where 
$$\lambda_{t}(a) = \sum_j t^j \lambda^j(a), \mbox{  } s_{t}(a) = \sum_j t^j S^j(a).$$
and $S^j$ denotes the $j$'th symmetric product.
Furthermore, we have for any $a,b \in
K_{\alg}^0(|{\mathcal M}_\sigma|^c)$ that 
\[\lambda_t(a-b)^{-1} = \lambda_t(b) s_{-t}(a),\] 
If $\rank ({\mathcal N}_i) = r_i$, 
we
have that 
\[\lambda_{-1}(\rank({\mathcal N}))^{-1} = \prod_{i=1}^{m-1}
(1-x^i)^{-r_i}.\] 
Combining the above formulae, we get the following Proposition.

\begin{prop}\label{repformula}
As an element  in 
$\lambda_{-1}^c({\mathcal M}_\sigma)^{-1}\in K_{\alg}^0(|{\mathcal M}_\sigma|^c)\otimes R({\Bbb Z}_m)$ we have that

\begin{eqnarray} (\lambda_{-1}^c{\mathcal M}_\sigma)^{-1} & = & 
\prod_{i=1}^{m-1}
(1-x^i)^{-r_i}  \sum_{i_1 , j_1}
S^{i_1}\left( \sum_{j=1}^{m-1}\tilde{T}^*_c\otimes x^j
\right)\nonumber \\ 
& & (-1)^{j_1}\lambda^{j_1}\left( \sum_{j=0}^{m-1} \left\{ \frac{1}{m}
\sum_{s=1}^l\sum_{\nu=0}^{m-1} \tilde{E}_{p_s}^\nu\otimes (\mu_m^{-\nu} 
(n_s) - \mu_m^{j-\nu}(n_s))\right\} \otimes 
x^j\right)  \\  \nonumber
& =& \prod_{i=1}^{m-1}(1-x^i)^{-r_i}
        \sum_{i=0}^{m-1}\left\{  \mbox{ }\sum_{\stackrel{1 (i_1+j_1)
\ldots+(m-1)(i_{m-1}+j_{m-1})}{\mbox{ } = \mbox{ } i \mod m}} \mbox{ }\prod_{q=1}^{m-1} 
S^{i_q}(\tilde{T}^*_c)\right.\\ \nonumber
& & (-1)^{j_q} \lambda^{j_q}\left. \left( \frac{1}{m}
\sum_{s=1}^l\sum_{\nu=0}^{m-1} \tilde{E}_{p_s}^\nu\otimes (\mu_m^{-\nu} 
(n_s) - \mu_m^{q-\nu}(n_s)) \right) \right\} 
\otimes x^i. \label{scf} 
\end{eqnarray}

\end{prop}

By identifying the fibre of $E_{p_s}$ at any point in 
$|{\mathcal M}_\sigma|^c$ with $\mathfrak{g}_{\Bbb C}$, our constructed lift of $\Id\times f$ will induce an automorphism
of the form $\Ad(g) : \mathfrak{g}_{\Bbb C} \rightarrow \mathfrak{g}_{\Bbb C}$ for some 
$g\in c^{-k_s}_s$. Thus the eigenspace splitting of $E_{p_s}$ is identified with the 
eigenspace decomposition of $\mathfrak{g}_{\Bbb C}$ under $\Ad(g)$. 

Let us describe how $\rank (E_{p_s}^\nu)$ can be expressed in terms of $
c_1, \ldots c_l$.
Let $T$ be a maximal torus through $g$ in $G$. Let $R$ be the set of roots of 
$G$. We then have the usual decomposition
\[\mathfrak{g}_{\Bbb C}=\mathfrak{t}_{\Bbb C}\oplus \bigoplus_{\alpha\in R}\mathfrak{g}^
\alpha_{\Bbb C}.\]
Here $\mathfrak{t}_{\Bbb C}$ is the Cartan subalgebra corresponding to $T$ and 
$\mathfrak{g}^\alpha_{\Bbb C}$ is the root-subspace of $\mathfrak{g}_{\Bbb C}$ 
corresponding to $\alpha\in R$. With respect to this decomposition $g$ acts trivially on $\mathfrak{t}_{\Bbb C}$ and it  multiplies by $\expo(2\pi i \alpha(c^{-k_s}_s))$ on 
$\mathfrak{g}^\alpha_{\Bbb C}.$

If we now let 
\[R_{s}^i = \left\{\alpha\in R| \expo(2\pi i \alpha(c^{-k_s}_s)) = 
e^{2\pi i\frac{i}{m}}
\right\},\]
we get that
\[R = \bigcup_{i=0}^{m-1} R_{s}^i,\]
and $\rank (E^i_{p_s}) = 2 \# R^i_{s}$ for all $i\neq 0 \mbox{ mod} 
\ m$.
 Moreover,
$\rank (E^0_{p_s}) = 2 \# R^0_{s} + \rank G$. We shall use the following
notation $2\# R_{s}^i = r^i_s$.

Let us compute the ranks $r_i$ defined above.
Let $A$ be a flat connection in $P$ invariant under the explicit lift 
${\phi}$ of $f$ 
representing $\delta$, such that $[A]\in|{\mathcal M}_\sigma|^c$. We then have that the cyclic group 
$\langle f\rangle$  of order $m$ acts on $T_{[A]}{\mathcal M}_\sigma\cong 
H^{0,1}(\Sigma,\overline{\partial}_A),$ where $\overline{\partial}_A$ is the 
$\overline{\partial}$-operator on the complexified adjoint bundle and we have a 
unique decomposition
\[H^{0,1}(\Sigma,\overline{\partial}_A) = \bigoplus_{i=0}^{m-1} V_{i},\]
where $V_i$ is the $e^{2\pi i \frac{i}{m}}$-eigenspace of $f$.

We have that $r_i = \dim_{\Bbb C} V_i$. As we have seen before,
\[r_i = \frac{1}{m} \sum_{\beta=0}^{m-1} \tr(f^\beta) e^{-2\pi i\frac{i\beta}{m}
}.\]
By the Lefschetz fixed point Theorem, we have
\begin{eqnarray*}\tr(f^\beta) &=& \sum_{s=1}^l \tr({\phi}^\beta(p_s)) 
(1-\expo(2\pi i\frac{n_{s}\beta}{m}))^{-1}\\
& =& \sum_{s=1}^l \left( \rank G + 2\sum_{\alpha\in R}\expo(2\pi i \beta 
\alpha(c^{-k_s}_s))\right) (1-\expo(2\pi i\frac{n_{s}\beta}{m}))^{-1}
\end{eqnarray*}
for $\beta\neq 0 \mod m$ and 
\[\tr(f^0) = \dim G (g-1).\]
So we get that
\begin{eqnarray*}
r_i & =& \frac{1}{m}\left( -\sum_{s=1}^l \sum_{\beta=1}^{m-1} \left( \rank G + 
2\sum_{\alpha\in R}\expo(2\pi i \beta \alpha(c^{-k_s}_s))\right)\frac{ e^{-2\pi i
\frac{i\beta}{m}}}{ (1-\expo(2\pi i
\frac{n_s\beta}{m}))} + \dim G (g-1)\right) \\
& =& \frac{1}{m}\left( \sum_{s=1}^l \left( \mu_m^i(n_s)\rank G  +
\sum_{j=1}^{m-1}r_s^j \mu_m^{i-j}(n_s)\right) + \dim G(g-1) 
\right) .\end{eqnarray*}

Implementing the formula we now have for $(\lambda_{-1}^c {\mathcal M}_\sigma)^{-1}$ 
into (\ref{csc}) we get 

\begin{theorem}\label{SAET2}
We have the following formula for the perturbative contribution from
 ${\mathcal M}(\Sigma_f)^c$ (when it is smooth)
\begin{equation}
P_c(k) = \frac{k^{-d_c}}{|Z_\delta|}\prod_{i=1}^{m-1}
(1-e^{2\pi i\frac{i}{m}})^{-r_i} 
  \expo(km\omega^{c}) \cup \lambda_c^{-1} \cup \Td(T_c)  \cap [{\mathcal
 M}(\tilde{\Sigma}',c(\delta))]. \label{SF}
\end{equation}
where
\begin{eqnarray*}
\lambda_c^{-1} = && \sum_{i=0}^{m-1}\left\{  
\sum_{ \stackrel{1 (i_1+j_1)
\ldots+(m-1)(i_{m-1}+j_{m-1})}{\mbox{ } = \mbox{ } i \mod m}} \mbox{ }
\bigcup_{q=1}^{m-1} (-1)^{j_q}\Ch^\cdot(S^{i_q}(\tilde{T}^*_c))\right.\cup\\
&&
\mbox{ } \Ch^\cdot \lambda^{j_q} \left. \left(  \frac{1}{m}
\sum_{s=1}^l\sum_{\nu=0}^{m-1} \tilde{E}_{p_s}^\nu\otimes (\mu_m^\nu(n_s
) - \mu_m^{q-\nu}(n_s)) \right)  \right\} e^{2\pi i\frac{i}{m}}
\end{eqnarray*}

and
$$
r_i= \frac{1}{m}\left( \sum_{s=1}^l \left( 
\mu_m^i(n_s)\rank G  +
\sum_{j=1}^{m-1}r_s^j \mu_m^{i-j}(n_s)\right) + \dim G(g-1) 
\right).
$$
Hence we see that the total contribution from the smooth component
${\mathcal M}(\Sigma_f)^c$ to the invariant of $\Sigma_f$ only depends on $m$, the
Seifert invariants of $\Sigma_f$ and $c$.
\end{theorem}

\begin{remark} In \cite{AH} we use this formula to get an expression for the contribution from a smooth component ${\mathcal M}(\Sigma_f)^c$ to the leading order term of the perturbation expansion of $Z^{(k)}_G(\Sigma_f)$. We see that the level shift $k \mapsto k+ h$ in the leading order term comes from the term $(\lambda_{-1}^c {\mathcal M}_\sigma)^{-1}$ and from the relation 
$$ \Td(T_c) = e^{\frac12 c_1(|{\mathcal M}_\sigma|^c) }\hat{A}(|{\mathcal M}|^c).$$
\end{remark}

\end{document}